\documentclass[preprint,12pt]{elsarticle}

\usepackage{mathtools}
\usepackage{graphicx}
\usepackage{amsfonts}
\usepackage{amssymb}
\usepackage{amsthm}
\usepackage{easyReview}
\usepackage[numbers]{natbib}
\usepackage{xargs}
\usepackage{color}
\usepackage[hidelinks]{hyperref}
\usepackage{xurl}
\usepackage{enumitem}
\usepackage{pgffor}

\theoremstyle{plain}
\newtheorem{teo}{Theorem}
\newtheorem{prop}[teo]{Proposition}

\newtheorem{lema}[teo]{Lemma}

\theoremstyle{definition}
\newtheorem{defi}[teo]{Definition}

\newtheorem*{teo*}{Theorem}

\theoremstyle{remark}
\newtheorem{remk}{Remark}

\newcommand{\pinterno}[2]{\left\langle #1 , #2 \right\rangle}
\newcommand{\normaSobolev}[2]{\left \| #1 \right\|_{W^{1,2}(#2)}}

\newcommand{\bb}[1]{\mathbb{#1}}
\newcommandx{\espHolder}[3][1=$k$, 2=$\lambda$, 3=$\Omega$, usedefault]{C^{#1,#2}(\overline{#3})}

\newcommand{\rquad}[1]{\foreach \n in {1,...,#1}{\quad}}
\newcommand{\negarquad}[1]{\foreach \n in {1,...,#1}{\mkern-18mu}}

\newcommand{\mytag}[2]{
	\text{#1}
	\@bsphack
	\begingroup
	\@onelevel@sanitize\@currentlabelname
	\edef\@currentlabelname{
		\expandafter\strip@period\@currentlabelname\relax.\relax\@@@
	}
	\protected@write\@auxout{}{
		\string\newlabel{#2}{
			{#1}
			{\thepage}
			{\@currentlabelname}
			{\@currentHref}{}
		}
	}
	\endgroup
	\@esphack
}

\journal{ArXiv }

\begin{document}

\begin{frontmatter}



\title{Homoclinic Solution to Zero of a Non-autonomous, Nonlinear, Second Order Differential Equation with Quadratic Growth on the Derivative}


\author[luiz]{Luiz Fernando de Oliveira Faria\fnref{label2}}
\fntext[label2]{Partially financed by FAPEMIG  APQ-02374-17, APQ-04528-22, RED-00133-21,  and   CNPq}
\affiliation[luiz]{organization={Departamento de Matemática, Universidade Federal de Juiz de Fora},
            city={Juiz de Fora},
          	postcode={CEP: 36036-900},
            state={MG},
            country={Brazil}}

\author[pablo]{Pablo dos Santos Corrêa Junior\corref{cor1}}
\ead{p.correa@usp.br}
\cortext[cor1]{Financed by Universidade Federal de Juiz de Fora, through scholarship}
\affiliation[pablo]{organization={Departamento de Matemática, Universidade de São Paulo},
	city={São Carlos},
	postcode={CEP: 13566-590},
	state={SP},
	country={Brazil}}

\begin{abstract}
This work aims to obtain a positive, smooth, even, and homoclinic to zero (i.e. zero at infinity) solution to a non-autonomous, second-order, nonlinear differential equation involving quadratic growth on the derivative. We apply Galerkin's method combined with Strauss' approximation on the term involving the first derivative to obtain weak solutions. We also study the regularity of the solutions and the dependence on their existence with a parameter. 
\end{abstract}




\end{frontmatter}


	\section{Introduction}
\label{intro}

The existence of homoclinic solutions for autonomous and non-autonomous differential equations and Hamiltonian systems is a crucial subject in qualitative theory (see \cite{ZHU20081975}).

In this work, the second-order equation in the real line is considered
\begin{equation}\label{problema-principal-falandoSobre}
	\begin{cases}
		-(A(u)u')'(t) + u(t)= \lambda a_1(t) |u(t)|^{q-1} +  |u(t)|^{p-1} + g(|u'(t)|),\quad \text{in}\,\, \bb{R} \\
  \;\;\;u(t)> 0,\;\;\;\rquad{22} \text{in}\,\, \bb{R} \\
		\;\;\;\lim_{t \rightarrow \pm \infty} u(t)= 0,
	\end{cases}
\end{equation}
with
\begin{description}
	\item[$(H_1)$] $1<q<2<p<+\infty$ and $a_1\in L^ \mathfrak{s}(\bb{R})\cap C(\bb{R})$, $\mathfrak{s}=\frac{2}{2-q}$ a positive even function;
	\item[$(H_2)$] $A:\bb{R}\rightarrow\bb{R}$ a Lipschitz, smooth (at least $C^1(\bb{R})$), non-decreasing function satisfying:
	\[\exists\;\gamma\in (0,1)\mbox{ such that }0<\gamma\leq A(t)\quad\forall t\in \bb{R};\]
	\item[$(H_3)$] $g:\bb{R}\rightarrow \bb{R}$ a continuous function satisfying:
	\begin{equation}\label{properties of g}
		0\leq sg(s)\leq |s|^\theta \mbox{ for all }s\in \bb{R},\mbox{ where }2<\theta\leq 3.
	\end{equation}
\end{description}

The equation in \eqref{problema-principal-falandoSobre} arises in several real phenomena, for instance, as the study of traveling wavefronts for parabolic reaction-diffusion equations with a local reaction term,  chemical models, and others, as mentioned in \cite{2005_Malaguti,2007_Marcelli}, and generalizes several classical equations such as Duffing-type equations \cite{kazarinoff1984nonlinear,
salas2014exact, 2011_Alves} or Linard-like systems \cite{zhang2012homoclinic}. 



Now, we state our main result.
\begin{teo}[Main problem]
	\label{teorema-central}
	There exists $\lambda^*>0$ such that, for all $\lambda\in (0,\lambda^*]$,  problem (\ref{problema-principal-falandoSobre}) has an even, positive and $C^2(\bb{R})$ homoclinic solution to the origin.	Moreover, as $\lambda\rightarrow 0$, this solution goes to $0$ in $C^0(\bb{R})$.
\end{teo}

We also find an additional result with respect to appropriate ranges of $\lambda$ in
order to guarantee the existence of  solutions.

\begin{prop} \label{prop1}
Assume the hypotheses of Theorem \ref{teorema-central}. If  $\lambda>0$ is sufficiently large, then equation \eqref{problema-principal-falandoSobre} has no (positive) solution in $ H^1(\mathbb{R})$. 
\end{prop}

The idea to consider problem \eqref{problema-principal-falandoSobre} came from article \cite{2011_Alves}, where the authors considered  a similar equation but with a different set of hypotheses; namely, their formulation was focused on the study of the equation
\[\begin{cases}
	-(A(u)u')' + u(t) = h(t,u(t)) + g(t,u'(t))\mbox{ in } \bb{R}\\
	u(\pm \infty)= u'(\pm \infty)=0,
\end{cases}\]
with
\begin{description}
	\item[($\widetilde{H}_1$)] $h,g:\bb{R}^2\rightarrow \bb{R}$ locally Hölder continuous, even on the first variable and $h(t,0)=g(t,0)=0$;
	\item[($\widetilde{H}_2$)] there exist constants $0<r_1,r_2<1$ and smooth functions $b\in L^1(\bb{R})\cap L^\infty(\bb{R})$ with $b(t)>0$ for all $t\in \bb{R}$, $a_1\in L^2(\bb{R})$ and $a_2\in L^{\frac{2}{1-r_2}}(\bb{R})$, satisfying
	\[b(t)|\mu|^{r_1}\leq h(t,\mu)\leq a_1(t)+a_2(t)|\mu|^{r_2},\quad \forall (t,\mu) \in \bb{R}^2;\]
	\item[($\widetilde{H}_3$)] there exist a constant $0<r_3<1$ and smooth functions $a_3\in L^{\frac{2}{1-r_3}}(\bb{R})$ and $a_4\in L^2(\bb{R})$ satisfying
	\[0\leq g(t,\eta)\leq a_4(t) + a_3(t)|\eta|^{r_3}\quad \forall (t,\eta)\in \bb{R}^2;\]
	\item[($\widetilde{H}_4$)] the function $A$ is smooth, nondecreasing and there exists $\gamma\in (0,1)$ satisfying
	\[0<\gamma\leq A(t)\quad \forall t\in \bb{R}.\]
\end{description}

By comparing with our work, we considered sup-linear growth on $u$ and $u'$, terms involving this type of growth are not covered in \cite{2011_Alves}. Another aspect that we would like to emphasize is the weakening of the hypothesis over $g$: comparing with \cite{2011_Alves}, we asked only for continuity over $g$, instead of Hölder continuity. Although the formulation presented here is not an immediate consequence of \cite{2011_Alves}, some techniques therein proved to be quite solid and very useful in the study of this type of problem, transcending the circumstances framed by the authors.

Our formulation presented some interesting challenges, for instance, the problem is not variational. Among the non-variational techniques, we chose the Galerkin Method as a tool to gather information about the existence of weak solutions. Although proving itself beneficial, the Galerkin Method presented us with other types of challenges to circumvent. For example, the nonlinear term $g(|u'(t)|)$ with $0\leq sg(s)\leq |s|^{\theta}$ and $2<\theta\leq 3$ enables us to take $g(s)\equiv sign(s)|s|^{2}$. Thus estimations involving $\int_\Omega|u'|^2$ become essential to the calculations but, at the same time, we cannot say much about it \emph{a priori}: this is due to the lack of information about $u'$, since the embedding theorems of $H_0^1(\Omega)$ do not provide substation information about $u'$ as they do for $u$. 

We consider the case $\theta = 3$ as the \emph{critical} one and treat it separately in our estimations. For $\theta>3$ we would get expressions involving $\int|u'|^{\theta-1}$ that we could not control, because $\theta-1>2$ and we only know that $u'\in L^2$; for this reason, we limited $\theta\leq 3$, and $\theta>2$ was required because we wanted to focus on the sup-linear case.

There is some literature about equations on domains in $\bb{R}^n$ involving the term $|\nabla u|^2$ in the nonlinearity (see \cite{2006_Abdellaoui,2007_Abdellaoui,2009_Pedro-gradiente,2013_Jeanjean}), some authors call this type of growth: ``\emph{critical growth on the gradient}''. Simple changes in how this term appears in the equation can have dramatic effects on the outcome. For instance, a simple change in the sign of $|\nabla u|^2$ can lead to a total failure to obtain a solution (even in the weak sense), see the article \cite{2007_Abdellaoui} for more information. We also would like to emphasize article \cite{2013_Jeanjean} for its results and broader discussion about PDE with quadratic growth on the gradient: the model problem studied by the authors is
\[-div(A(x)\nabla u)=c_0(x)u+ \mu(x)|\nabla u|^2 +f(x),\]
with suitable hypothesis. In this context, our problem (\ref{problema-principal-falandoSobre}) presents a similar structure that was not covered before, thus -we believe- it contributes to the discussion previously mentioned.

The methods applied in our work require certain symmetry, which is due mainly to a lack of a comparison principle (known to the authors) to guarantee that some limit-functions are not zero almost everywhere (a.e.), (see Proposition \ref{proposicao-propiedades-un}). To overcome this obstacle, we founded this work focusing on the set $\bb{E}^1_0(I)=\{u\in H^1_0(I);u(t)=u(-t)\mbox{ a.e.}\}$, $I\subset \bb{R}$ a symmetric interval, which is the subset of $H^1_0(I)$ consisting of even functions. This set can be understood as the set of radial symmetric functions in $\bb{R}$.

In order to develop our study, in Section \ref{sec:Solution in a bounded interval}, we started by analyzing our equation on a bounded interval:
{
	\begin{equation}\label{main problem-INTRODUCAO}
		\begin{cases}
			-(A(u)u')'(t) + u(t)= \lambda a_1(t) |u(t)|^{q-1} +  |u(t)|^{p-1} + g(|u'(t)|),\mbox{ in } (-n,n) \\
			u( n)=u(- n) =0.
		\end{cases}\tag{$P_n$}
	\end{equation}
}
This restriction was essential to realize our estimations and to obtain upper bound constants that were crucial to construct the solution, in $\bb{R}$, to the problem (\ref{problema-principal-falandoSobre}). The process developed in Section \ref{sec:Solution in a bounded interval} consisted mainly of two steps:
\begin{description}
	\item[\textbf{First}] Approximate $g$ by a sequence of Lipschitz functions $(f_k)$ using the \emph{Strauss Approximation}; this sequence received this name after its first  appearance in the article \cite{1970_Strauss}.

		This approximation was useful because it helped us to work with the necessary estimations without extra hypotheses over $g$.  We followed \cite{2020_Araujo} in the definition and presentation of the properties of the sequence $(f_k)$. In this article, the authors used this approximation to avoid the usage of the Ambrosetti-Rabinowitz condition and were able to obtain a positive solution to the equation
	\[\begin{cases}
		-\Delta u = \lambda u^{q(r)-1} + f(r,u)&\mbox{ in } B(0,1)\\
		u>0 \quad &\mbox{in } B(0,1)\\
		u=0 \quad &\mbox{on } \partial B,
	\end{cases}\]
	see \cite{2020_Araujo} for more information.

	We would like to emphasize that, in \cite{2020_Araujo}, the authors used this approximation in a term involving $u$; namely they used it to approximate $f(r,u)$. In our work, we used it in $u'$. 
 
	\item[\textbf{Second}] We used the sequence $(f_k)$ to define an approximate problem in $(-n,n)$ and used the Galerkin Method to obtain a weak solution. Then, using the work done by Gary M. Liberman \cite{1988_Lieberman}, we obtained the \emph{a priori estimation} summarized in Proposition \ref{proposicao-regularidade-solucao-problema-aproximado}. Thus we obtained a strong solution to this problem. Afterward, we were able to construct a strong solution to the problem (\ref{main problem-INTRODUCAO}).
\end{description}

In Section \ref{secao: solucao na reta toda} we used the pieces of information gathered in Section \ref{sec:Solution in a bounded interval} to construct a solution to the problem (\ref{problema-principal-falandoSobre}), thus proving Theorem \ref{teorema-central}. We also prove Proposition \ref{prop1} in the Section \ref{secao: solucao na reta toda}.  We would like to point out the role of Section \ref{sec: asympotic solution to problem Pn}: there we study the asymptotic behavior, in respect to $\lambda$, of the solution to the problem (\ref{main problem}) -- the arguments presented were inspired by the article \cite{2021_Faria}. This was useful to tackle the last assertion of Theorem \ref{teorema-central}.

\section{Solution in a bounded interval}
\label{sec:Solution in a bounded interval}

First, we will obtain a solution to a problem related to (\ref{problema-principal-falandoSobre}); namely, we will study
{ \begin{equation}\label{main problem}
		\begin{cases}
			-(A(u)u')'(t) + u(t)= \lambda a_1(t) |u(t)|^{q-1} +  |u(t)|^{p-1} + g(|u'(t)|),\mbox{in } (-n,n) \\
			u( n)=u(- n) =0,
		\end{cases}\tag{$P_n$}
\end{equation}}
with the same set of hypothesis $(H_1),(H_2)$ and $(H_3)$. The motivation for this approach is to construct a solution to the problem (\ref{problema-principal-falandoSobre}) using the solutions of (\ref{main problem}). Although the analysis of (\ref{main problem}) is easier since it is defined over $(-n,n)$, rather than $\bb{R}$, the lack of hypothesis over $g$ creates a difficult situation for our estimations. To overcome this matter, we will utilize the \emph{Strauss Approximation} on $g$ at the same time that we approximate the problem (\ref{main problem}). Let us define the sequence of functions that will approximate $g$.

Define $G(s)= \int_0^s g(t)dt$ so that $G$ is differentiable and $G'(s) = g(s)$. By means of $G$ we shall construct a sequence of approximations of $g$ by Lipschitz functions $f_k:\mathbb{R}\longrightarrow \mathbb{R}.$
Let
\begin{equation}
	f_k(s) = \begin{cases}
		-k[G(-k-\frac{1}{k}) - G(-k)], &\text{if} \,\, s\leq -k\\
		-k[G(s-\frac{1}{k}) - G(s)],  &\text{if}\,\, -k\leq s \leq \frac{-1}{k}\\
		k^2s[G(\frac{-2}{k}) - G(\frac{-1}{k})],  &\text{if}\,\,\, \frac{-1}{k}\leq s \leq 0 \\
		k^2s[G(\frac{2}{k}) - G(\frac{1}{k})],  &\text{if}\,\, 0\leq s \leq \frac{1}{k} \\
		k[G(s+\frac{1}{k}) - G(s)],  &\text{if}\,\, \frac{1}{k}\leq s \leq k \\
		k[G(k+\frac{1}{k}) - G(k)], &\text{if} \,\, s\geq k. \\
	\end{cases}
\end{equation}
\begin{remk}
	The construction of the sequence $(f_k)$ is due to \cite{1970_Strauss}.
\end{remk}
The advantage of this sequence lies in the properties that one can obtain from it:
\begin{teo}[{\cite[Lemma~1]{2020_Araujo}}]\label{teo of sequence f_k}
	The sequence $f_k$ as defined above satisfies:
	\begin{enumerate}
		\item $sf_k(s)\geq 0$ for all $s \in \mathbb{R}$;
		\item for all $k \in \mathbb{N}$ there is a constant $c(k)$ such that $|f_k(\xi) - f_k(\eta)| \leq c(k) |\xi -\eta|$, for all $\xi,\eta \in \mathbb{R}$;
		\item $f_k$ converges uniformly to $g$ in bounded sets.
	\end{enumerate}
\end{teo}
\begin{remk}
From the definition of the sequences $f_k$, and the fact that $sign(g(s))$ = $sign(s)$ for all $s \in \mathbb{R}$,  it follows without difficulties that \textit{1} is true. In \cite[Pag.~6, Prop.~5]{2020_Araujo} one can find a detailed proof of \textit{2}, so we will only prove \textit{3} by an alternative argumentation.
\end{remk}
\begin{proof}
	Given a bounded set $J\subset \mathbb{R}$, there exists $m_0\in \mathbb{N}$ such that $J\subset(-m_0,m_0)$; so to prove \textit{3} we only need to prove that it holds in intervals such as $(-m,m)$, $m\in \mathbb{N}$. We may also assume that $k>m$. Consider the following cases:

	\textbf{Case I.} $ -m<s\leq \frac{-1}{k}$
	Here we have that
	\begin{equation*}
		|f_k(s) - g(s)| = \left|-k\left[G\left(s-\frac{1}{k}\right) - G(s)\right] -g(s)\right| = \left|\frac{\left[G(s-\frac{1}{k}) - G(s)\right]}{\frac{-1}{k}} - g(s)\right|.
	\end{equation*}
	Then, by the Fundamental Theorem of Calculus, we conclude that $f_k\rightarrow g$ uniformly.

	\textbf{Case II.} $\frac{-1}{k}\leq s\leq 0$

	Since $g(0)=0$ and $g$ is continuous, given $\epsilon>0$ there exists $\delta >0$ such that $|t|<\delta$ implies $|g(t)|<\epsilon/2$. Let $k_0\in \mathbb{N}$ be such that $k_0>m$ and $k_0>2/\delta$. Then, for $k>k_0$
	\begin{align*}
		|f_k(s) - g(s)| &= \left|k^2s\left[G\left(\frac{-2}{k}\right) - G\left(\frac{-1}{k}\right)\right] -g(s) \right|\\
		&\leq k^2|s|\left|\int_{-1/k}^{-2/k}|g(t)|dt\right| + |g(s)| \\
		&\leq k^2 (\frac{1}{k^2}) \frac{\epsilon}{2} + \frac{\epsilon}{2} = \epsilon \quad \forall s \in [\frac{-1}{k},0].
	\end{align*}
	Proving the desired convergence.

	For the cases $0\leq s \leq \frac{1}{k}$ and $\frac{1}{k}\leq s< m$ the arguments are similar.

\end{proof}
\begin{lema}[{\cite[Lemma~2]{2020_Araujo}}]\label{upper-estimative for the seq f_k}
	Let $g:\mathbb{R}\rightarrow\mathbb{R}$ be a continuous function satisfying (\ref{properties of g}). Then the sequence $f_k$ of Theorem \ref{teo of sequence f_k} satisfies
	\begin{enumerate}
		\item For all $ k \in \mathbb{N}$, $0\leq sf_k(s)\leq C_1 |s|^\theta$ for every $|s|\geq \frac{1}{k}$;
		\item for all $ k \in \mathbb{N}$, $0\leq sf_k(s)\leq C_1 |s|^2$ for every $|s|\leq \frac{1}{k}$;
	\end{enumerate}
	where $C_1$ is a constant independent of $k$.
\end{lema}
\begin{proof}

	See \cite[Pag.~8, Lemma~2]{2020_Araujo}.
\end{proof}

Now, we are in condition -- using $(f_k)$ -- to define a problem that approximates problem (\ref{main problem}). Let $\psi \in L^2(-n,n)$ be a positive, even function. We define our \emph{approximate problem} by:
{ \begin{equation}\label{problem with f_k}
		\begin{cases}
			-(A(u)u')'(t) +u(t) = \lambda a_1(t) |u(t)|^{q-1} +  |u(t)|^{p-1}  + f_k(|u'(t)|) + \frac{\psi}{k},\mbox{in} {\small (-n,n)} \\
			u(n) = u(-n) = 0.
		\end{cases}\tag{$P^k_n$}
\end{equation}}
In the next subsection, we will utilize the \emph{Galerkin Method} to obtain a solution to (\ref{problem with f_k}); afterward, we will let $k$ vary and thus, as $k\rightarrow \infty$, obtain a solution to (\ref{main problem}). Before jumping into the next subsection, let us define what we understand as \emph{weak solution} to problem (\ref{main problem}):
\begin{defi}
	We will call $w\in H_0^1(-n,n)$ a \emph{weak solution} of (\ref{main problem}) if
	\[
	\int_{-n}^{n}A(w)w'v' + \int_{-n}^{n}wv = \int_{-n}^{n}\lambda a_1|w|^{q-1}v + \int_{-n}^{n}|w|^{p-1}v +\int_{-n}^{n}g(|w'|)v
	\]
	for all $v\in H_0^1(-n,n).$
\end{defi}
\begin{remk}
	We will use, for the sake of clarity, the notation $\|\cdot\|_{W^{1,2}}$ for the usual norm of $H_0^1$ and for ($\|u\|_{L^2}+\|u'\|_{L^2}$) or $(\|u\|^2_{L^2}+\|u'\|^2_{L^2})^{1/2}$. Since these norms are equivalent the results will not change but the constants may. In most of the cases $\|u\|_{W^{1,2}}=\|u\|_{L^2}+\|u'\|_{L^2}$ . We also emphasize that, when the context is clear, we will omit the domain in norms such as those from the spaces $L^p(-n,n)$.
\end{remk}
\begin{remk}
	The integrals in the definition above are well defined, see for instance the estimations of Proposition \ref{proposicao-solucao-fraca-problema-aproximado}. The same is true for the definitions given in the next subsection.
\end{remk}

\subsection{Solution to the approximate problem}
\label{subsec:Solution to the approximate problem}

Our main goal in this subsection is to prove the following theorem:

\begin{teo}\label{teo-da-solucao-do-problema-aproximado}
	There exist $\lambda^*>0$, $\beta\in (0,1)$ and $k^*\in \mathbb{N}$ for which the problem (\ref{problem with f_k}) admits a nontrivial, even, non-negative $C^{1,\beta}[-n,n]\cap C^2(-n,n)$ solution  for every $\lambda\in (0,\lambda^*)$ and $k\geq k^*$.
\end{teo}
As mentioned, we will utilize the \emph{Galerkin Method}; thus we will start by presenting the foundations that this method requires. The next lemma is a well-known result, but it plays a central role in all arguments involving the Galerkin Method.
\begin{lema}\label{corl fixed point brower}
	Let $\mathfrak{F}:\mathbb{R}^N\rightarrow\mathbb{R}^N$ be a continuous function such that $\langle \mathfrak{F}(x),x\rangle\geq 0$ for all $x\in \mathbb{R}^N$ with $\|x\|_{\mathbb{R}^N} =r$. Then there exists $x_0$ in the closed ball $B[0,r]$ such that $\mathfrak{F}(x_0)=0$.
\end{lema}
\begin{proof}
	See \cite[Chap.~5, Theorem~5.2.5]{1989_Kesavan_BOOK}.
\end{proof}
Now we will define an entity called \emph{E-weak solution}. It is well known that the main focus of the Galerkin Method is to obtain a weak solution, but we will utilize it to obtain an E-weak solution first.
\begin{defi}
	A function $w\in H^1_0(-n,n)$ is called an \emph{E-weak solution} of (\ref{problem with f_k}) if $w$ is an \emph{even function} satisfying
	\begin{multline*}
		\int_{-n}^{n}A(w)w'\varphi' + \int_{-n}^{n}w\varphi= \int_{-n}^{n}\lambda a_1 |w|^{q-1}\varphi + \int_{-n}^{n}|w|^{p-1}\varphi\\ + \int_{-n}^{n}f_k(|w'|)\varphi + \int_{-n}^{n}\frac{\psi}{k}\varphi\end{multline*}
	for all $\varphi\in \bb{E}^1_0(-n,n)=\{u\in H^1_0(-n,n);u(t)=u(-t)\mbox{ a.e.}\}$.
\end{defi}
The use of the E-weak solution will be central in our argumentation to obtain an \emph{even} solution to the problem (\ref{problem with f_k}). This symmetry -- being even -- will also be beneficial in the use of our comparison principle, which is stated as follows:
\begin{teo}[{\cite[Theorem~3.1]{2011_Alves}}]\label{criterio-de-comparação}
	Let $\sigma:(0,+\infty)\rightarrow \bb{R}$ be a continuous function such that the mapping $(0,+\infty)\ni s \mapsto \frac{\sigma(s)}{s}$ is strictly decreasing and $\rho>0$. Suppose that there exist \emph{even functions} $v,w\in C^2(-\rho,\rho)\cap C[-\rho,\rho]$ such that:
	\begin{enumerate}
		\item $(A(w)w')'-w+\sigma(w)\leq 0 \leq (A(v)v')-v+\sigma(v)$ in $(-\rho,\rho)$;
		\item $v,w\geq 0$ in $(-\rho,\rho)$ and $v(\rho)\leq w(\rho)$;
		\item $\{x\in (-\rho,\rho);v(x)=0\}$ and $\{x\in(-\rho,\rho);w(x)=0\}$ have null measure in $\bb{R}$;
		\item $v'\cdot w'\geq 0$ in $(-\rho,\rho)$;
		\item $v',w'\in L^{\infty}(-\rho,\rho)$.
	\end{enumerate}
	Then $v\leq w$ in $(-\rho,\rho)$.
\end{teo}
\begin{proof}

	See \cite[Pag.~2419, Thm~3.1]{2011_Alves}.
\end{proof}

Turns out that, obtaining an E-weak solution enables us to recuperate a weak solution in the usual sense:
\begin{defi}
	\label{definicao-sol-fraca-problema-com-fk}
	We will call $w\in H_0^1(-n,n)$ a \emph{weak solution} of (\ref{problem with f_k}) if
	\begin{multline*}
		\int_{-n}^{n}A(w)w'v' + \int_{-n}^{n}wv = \int_{-n}^{n}\lambda a_1|w|^{q-1}v + \int_{-n}^{n}|w|^{p-1}v \\ +\int_{-n}^{n}f_k(|w'|)v +\int_{-n}^{n}\frac{\psi}{k}v
	\end{multline*}
	for all $v\in H_0^1(-n,n).$
\end{defi}

This is achieved by
\begin{lema}[{\cite[Lemma~4.1]{2011_Alves}}]\label{E-weak solution is also weak solution}
	Let $w\in H^1_0(-n,n)$ be an E-weak solution of (\ref{problem with f_k}). Then $w$ is a weak solution of (\ref{problem with f_k}).
\end{lema}
\begin{proof}

	See \cite[Pag.~2421, Lemma~ 4.1]{2011_Alves}.
\end{proof}

The set $\bb{E}_0^1(-n,n)$ can be understood as the set of radial symmetric functions in $\bb{R}$. One can prove without difficulties the following properties of $\bb{E}_0^1(-n,n)$:
\begin{enumerate}[label=\roman*)]
	\item it is a Hilbert Space;
	\item it is separable;
	\item it has an orthonormal basis.
\end{enumerate}
Let $\bb{E}^1_0(-n,n)=\{u\in H^1_0(-n,n);u(t)=u(-t)\mbox{ a.e.}\}$ and $(e_l)_{l=1}^{\infty}$ be an orthonormal   basis of $\bb{E}^1_0(-n,n)$.

Define $V_M =span \{e_1,\ldots,e_M\}$; then for every $u\in V_M$ there exist $\xi_1,\ldots,\xi_M $ in $\mathbb{R}$ such that $u = \sum_{i=1}^{M}\xi_ie_i$. By means of $T:V_M\rightarrow \mathbb{R}^M$, $T(u)= T(\sum_{i=1}^{M}\xi_ie_i) = (\xi_1,\ldots,\xi_M)$, which is a linear isomorphism and preserve norm, we may define $\mathfrak{F}:\mathbb{R}^M\rightarrow\mathbb{R}^M$ such that
\begin{equation}\label{function f in finite dimension}
	\mathfrak{F}(\xi) =(\mathfrak{F}_1(\xi),\ldots,\mathfrak{F}_M(\xi))
\end{equation}
and
\begin{equation*}
	\begin{array}{lll}
		\mathfrak{F}_j(\xi)&=& \int_{-n}^{n}A(u)u'e'_j + \int_{-n}^{n}ue_j - \int_{-n}^{n}\lambda a_1|u|^{q-1}e_j\\  &&- \int_{-n}^{n}|u|^{p-1}e_j -\int_{-n}^{n}f_k(|u'|)e_j -\int_{-n}^{n}\frac{\psi}{k}e_j,
	\end{array}
\end{equation*}
where $j\in \{1,\ldots,M\}$ and $u = T^{-1}(\xi)$, for all $\xi \in \bb{R}^M$.

\begin{lema}\label{Lemma f is continuous}
	The function $\mathfrak{F}$ is continuous.
\end{lema}
\begin{remk}
	Our proof will use the fact that, if $(x_n)$ is a sequence that converges to $x$ and, for all subsequence $(x_{n_l})$ of $(x_n)$, there exist a subsequence $(x_{n_{l_k}})$ of $(x_{n_l})$ such that $\mathfrak{F}(x_{n_{l_k}})$ converges to $\mathfrak{F}(x)$, then $\mathfrak{F}(x_n)$ converges to $\mathfrak{F}(x)$.
\end{remk}
\begin{proof}
	Given $\xi =(\xi_1,\ldots,\xi_M)\in \mathbb{R}^M$, let $(\xi_l)_{l=1}^{\infty}$ be a sequence in $\mathbb{R}^M$ such that $\|\xi_l-\xi\|_{\mathbb{R}^M}\rightarrow 0$. By means of $T$ we can identify $T^{-1}(\xi)=u=\sum_{i=1}^{M}e_i\xi_i$ and $T^{-1}(\xi_l)=u_l=\sum_{i=1}^{M}e_i\xi_i^l.$ Since $T$ is isometry we have that $\|u_l-u\|_{W^{1,2}}\rightarrow 0$. That is, $\|u_l-u\|_{L^2}\rightarrow 0$ and $\|u'_l-u'\|_{L^2}\rightarrow 0$. Taking a subsequence, if necessary, we may assume that
	\begin{align*}
		u_l(x) \rightarrow u(x) \;\text{a.e. on}\; (-n,n),\\
		u'_l(x)\rightarrow u'(x)\; \text{a.e. on}\; (-n,n),\\
	\end{align*}
	and $|u_l(x)|\leq h_1(x)$, $|u'_l(x)|\leq h_2(x)$ a.e. on $(-n,n)$, with $h_1,h_2\in L^2(-n,n)$. Let $j\in \{1,2,\ldots,M \}$, we will prove that $\mathfrak{F}_j(\xi_l)\rightarrow \mathfrak{F}_j(\xi)$.
	\begin{align}\label{f is continuous, first estimative}
		\left|\int_{-n}^{n}A(u_l)u'_le'_j - \int_{-n}^{n}A(u)u'e'_j\right|\leq \int_{-n}^{n}\left(|u'_l||A(u_l)-A(u)|+|A(u)||u'_l-u'|\right)|e'_j|,
	\end{align}
	since $|u'_l(x)||A(u_l(x))-A(u(x))||e'_j(x)|\rightarrow 0$ a.e. and $|A(u(x))||u'_l(x)-u'(x)||e'_j(x)|\rightarrow0$ a.e., by the Dominated Convergence Theorem (D.C.T) (\ref{f is continuous, first estimative}) tends to zero as $l\rightarrow +\infty$.
	\begin{equation}\label{f is continuous, second estimative}
		\left|\int_{-n}^{n}u_le_j-\int_{-n}^{n}ue_j\right|\leq \int_{-n}^{n}|u_l-u||e_j| \rightarrow 0 \quad\text{by (D.C.T)}.
	\end{equation}
	\begin{align}
		\left|\int_{-n}^{n}[\lambda a_1(|u_l|^{q-1}-|u|^{q-1}) + (|u_l|^{p-1}-|u|^{p-1}) + (f_k(u'_l)-f_k(u'))] e_j\right|\nonumber \\
		\leq \int_{-n}^{n}\lambda |a_1|\left||u_l|^{q-1}-|u|^{q-1}\right||e_j| + \int_{-n}^{n}\left||u_l|^{p-1}-|u|^{p-1}\right||e_j|\nonumber\\ +  \int_{-n}^{n}\left|f_k(|u'_l|)-f_k(|u'|)\right||e_j|,
	\end{align}
	since that $|u_l|^{q-1}\rightarrow |u|^{q-1}$ a.e. and $|u_l|^{p-1}\rightarrow |u|^{p-1}$ a.e., (D.C.T) implies that the first two integrals above converge to zero. Using the second item of Theorem \ref{teo of sequence f_k}, we have
	\begin{equation}\label{f is continuous, third estimative}
		\int_{-n}^{n}\left|f_k(|u'_l|)-f_k(|u'|)\right||e_j| \leq \int_{-n}^{n}c(k)|u'_l-u'||e_j|.
	\end{equation}
	Then, by (D.C.T), (\ref{f is continuous, third estimative}) converges to $0$ as $l\rightarrow +\infty$.

	These estimations show us that for every subsequence $(\xi_{l_k})$ of $(\xi_l)$, there exists a subsequence $(\xi_{l_{k_n}})$ of $(\xi_{l_k})$ that $\mathfrak{F}_j(\xi_{l_{k_n}})\rightarrow \mathfrak{F}_j(\xi)$. Therefore $\mathfrak{F}_j(\xi_l)\rightarrow \mathfrak{F}_j(\xi)$.
\end{proof}

\begin{prop}\label{proposicao-solucao-fraca-problema-aproximado}
	There exist $\lambda^{*}>0$ and $k^*\in \mathbb{N}$ for which the problem (\ref{problem with f_k}) admits a nontrivial weak solution for every $\lambda\in (0,\lambda^*)$ and $k\geq k^*$.
\end{prop}
\begin{remk}
	We will, in fact, search for an E-weak solution; but as seen in Lemma \ref{E-weak solution is also weak solution}, this will be also a weak solution.
\end{remk}
\begin{proof}
	Our aim is to use Lemma \ref{corl fixed point brower}, with the function $\mathfrak{F}$ defined in (\ref{function f in finite dimension}). Given $\xi \in \mathbb{R}^M$, we have that
	\begin{multline}\label{inerproduct fx}
		\langle \mathfrak{F}(\xi),\xi \rangle = \int_{-n}^{n}A(u)|u'|^2 + \int_{-n}^{n}|u|^2 - \int_{-n}^{n}\lambda a_1|u|^{q-1}u - \int_{-n}^{n}|u|^{p-1}u \\ -\int_{-n}^{n}f_k(|u'|)u -\int_{-n}^{n}\frac{\psi}{k}u.
	\end{multline}
	\phantom\qedhere
\end{proof}
In the following, we will estimate these integrals. We have that
\begin{align}
	\int_{-n}^{n}\lambda a_1|u|^{q-1}u &\leq \lambda \|a_1\|_{L^\mathfrak{s}(\mathbb{R})}\|u\|^q_{L^2} \leq \lambda C_2 \|u\|^{q}_{W^{1,2}},\\
	\int_{-n}^{n}\frac{\psi}{k}u &\leq \frac{\|\psi\|_{L^2(-n,n)}\|u\|_{W^{1,2}}}{k}.
\end{align}
Now let $\tilde{u}:\mathbb{R}\rightarrow \mathbb{R}$ be the extension by zero of $u$, then
\begin{align}
	\int_{-n}^{n}|u|^{p-1}u \leq \int_{-n}^{n}|u|^{p} &= \int_{-n}^{n}|u|^2|u|^{p-2}\\
	&\leq \|\tilde{u}\|_{L^{\infty}(\mathbb{R})}^{p-2}\int_{-n}^{n}|u|^{2}\\
	&=\|\tilde{u}\|_{L^{\infty}(\mathbb{R})}^{p-2}\|u\|^2_{L^2}\\
	&\leq C^{p-2}\|u\|^{p-2}_{W^{1,2}} \|u\|^2_{W^{1,2}} = C^{p-2}\|u\|^p_{W^{1,2}}.
\end{align}
Where  $C$ is the constant for the embedding $W^{1,2}(\mathbb{R})\hookrightarrow L^{\infty}(\mathbb{R})$.

Define
\[\Omega_> =\{s\in (-n,n); |u'(s)|\geq \frac{1}{k}\}\;\;\text{and}\;\;\Omega_<=\{s\in (-n,n); 0<|u'(s)|\leq \frac{1}{k} \}. \]
Then
\[\int_{-n}^{n}f_k(|u'|)u = \int_{\Omega_>}f_k(|u'|)u + \int_{\Omega_<}f_k(|u'|)u. \]
Notice that by Lemma \ref{upper-estimative for the seq f_k},
\begin{align*}
	\int_{\Omega_<}f_k(|u'|)u &\leq \int_{\Omega_<}C_1|u'||u|\leq \int_{\Omega_<}\frac{C_1}{k}|u|\\
	&\leq\frac{C_1}{k}\int_{-n}^{n} |u| \leq \frac{C_1(2n)^{1/2}}{k}\|u\|_{L^{2}}\\
	&\leq \frac{C_1(2n)^{1/2}}{k}\|u\|_{W^{1,2}}.
\end{align*}

To estimate the integral over $\Omega_>$, consider the following cases :
\newline
\textbf{Case 1.} $2<\theta<3$.

Using Lemma \ref{upper-estimative for the seq f_k}, we have
\begin{align*}
	\int_{\Omega_>}f_k(|u'|)u &\leq \int_{\Omega_>}C_1|u'|^{\theta-1}|u|\leq \int_{-n}^{n}C_1|u'|^{\theta-1}|u|\\
	&\leq C_1\left(\int_{-n}^{n}|u|^{w}\right)^{\frac{1}{w}}\left(\int_{-n}^{n}|u'|^2\right)^{\frac{\theta-1}{2}}\\
	&\leq C_1\left(\int_{\mathbb{R}}|\tilde{u}|^2|\tilde{u}|^{w-2}\right)^{\frac{1}{w}}\|u'\|^{\theta-1}_{L^2}\\
	&\leq C_1 \|\tilde{u}\|_{L^\infty(\mathbb{R})}^{\frac{w-2}{w}}\|u\|^{\frac{2}{w}}_{L^2}\|u'\|^{\theta-1}_{L^2}\\
	&\leq C_1C^{\frac{w-2}{w}}\|u\|^{\frac{w-2}{w}}_{W^{1,2}}\|u\|^{\frac{2}{w}}_{W^{1,2}}\|u\|^{\theta-1}_{W^{1,2}} = C_1C^{\frac{w-2}{w}}\|u\|^{\theta}_{W^{1,2}}.
\end{align*}
Where $w=\left(\frac{2}{\theta-1}\right)' = \frac{2}{3-\theta}>2$.
\newline
\textbf{Case 2.} $\theta=3$.
\begin{align*}
	\int_{\Omega_>}f_k(|u'|)u &\leq\int_{\Omega_>}C_1|u'|^2|u| \leq\int_{-n}^{n}C_1|u'|^2|u|\\
	&\leq C_1\|\tilde{u}\|_{L^{\infty}(\mathbb{R})}\|u'\|^2_{L^2} \leq C_1C\|u\|_{W^{1,2}}\|u\|_{W^{1,2}}^2\\
	&= C_1C\|u\|_{W^{1,2}}^3.
\end{align*}
Now we are able to estimate (\ref{inerproduct fx}). Notice that $\frac{w-2}{w}=\theta-2$.

\begin{align*}
	\langle \mathfrak{F}(\xi),\xi\rangle \geq &\gamma \|u\|^2_{W^{1,2}} -\lambda C_2 \|u\|^{q}_{W^{1,2}} -C^{p-2}\|u\|^p_{W^{1,2}}\\ &- C_1\max\{C^{\theta-2},C\}\|u\|^{\theta}_{W^{1,2}}  -\left(\frac{C_1(2n)^{1/2}}{k}+\frac{\|\psi\|_{L^2(-n,n)}}{k}\right)\|u\|_{W^{1,2}}.
\end{align*}
Define $Z_k:\mathbb{R}^+\rightarrow \mathbb{R}$ by
\begin{align*}
	Z_k(x) = \gamma x^2 - \lambda C_2 x^q-C^{p-2}x^p -C_1\max\{C^{\theta-2},C\}x^{\theta} -\left(\frac{C_1(2n)^{1/2}}{k}+\frac{\|\psi\|_{L^2(-n,n)}}{k} \right)x.  \end{align*}
We would like to find $r\in \mathbb{R}^+_*$ such that
\begin{equation}\label{x_1 that we want}
	\gamma r^2 -C^{p-2}r^p -C_1\max\{C^{\theta-2},C\}r^{\theta} > \frac{r^2}{2}\gamma
\end{equation}
or equivalently,
\[\frac{\gamma}{2} > C^{p-2}r^{p-2} +C_1\max\{C^{\theta-2},C\}r^{\theta-2} .  \]
For this, if we take
\[\delta_1 = \min\left\{\left(\frac{\gamma}{4C^{p-2}}\right)^{1/(p-2)},\left(\frac{\gamma}{4C_1\max\{C^{\theta-2},C\}}\right)^{1/(\theta-2)} \right\},\]
then for $0<r<\delta_1$ (\ref{x_1 that we want}) is true. Consequently,
\[Z_k(r)\geq \frac{r^2}{2}\gamma -\lambda C_2\delta_1^{q} -\left(\frac{C_1(2n)^{1/2}}{k}+\frac{\|\psi\|_{L^2(-n,n)}}{k}\right)\delta_1. \]
Define $\rho_1 = \frac{r^2}{2}\gamma -\lambda C_2\delta_1^{q}$. We will adjust $\lambda>0$ so that $\rho_1>0$; for this if $\rho_1>0$ it would imply that
\[\frac{r^2}{2}\gamma -\lambda C_2\delta_1^{q}>0 \Leftrightarrow \frac{r^2\gamma}{2C_2\delta^q_1}> \lambda. \]
Take $\Lambda^* = \frac{r^2\gamma}{2C_2\delta^q_1} $ and $0<\lambda<\Lambda^*$. Thus, $\rho_1>0$ and we can find $k^*\in \mathbb{N}$ such that for $k>k^*$, $\rho_1>\left(\frac{C_1(2n)^{1/2}}{k}+\frac{\|\psi \|_{L^2(-n,n)}}{k}\right)\delta_1>0$. Therefore, for $0<r<\delta_1$, $0<\lambda<\Lambda^*$ and $k>k^*$
\[ Z_k(r)>0, \]
and so, with $\|u\|_{W^{1,2}}= r$,
\begin{equation}\label{inner_product limitation, case theta<3}
	\langle \mathfrak{F}(\xi),\xi\rangle >0.
\end{equation}

By Lemma \ref{corl fixed point brower}, there exists $y_M\in B[0,r]$ such that $\mathfrak{F}(y_M) =0$ that is, identifying $v_M = T^{-1}(y_M)$, for all $j\in \{1,\ldots,M\}$

\begin{align}\label{weak solution depending on M}
	&\int_{-n}^{n}A(v_M)v'_Me'_j + \int_{-n}^{n}v_Me_j =\\ \nonumber
	&=\int_{-n}^{n}\lambda a_1|v_M|^{q-1}e_j + \int_{-n}^{n}|v_M|^{p-1}e_j +\int_{-n}^{n}f_k(|v'_M|)e_j + \int_{-n}^{n}\frac{\psi}{k}e_j.
\end{align}
Therefore (\ref{weak solution depending on M}) holds for all $\varphi\in V_M$, because $\{e_1,\ldots,e_M\}$ is a basis of $V_M$. Notice that
\begin{equation}
	\|v_M\|_{W^{1,2}}\leq r \mbox{ for all } M \in \bb{N}.
\end{equation}
\begin{remk}\label{observacao r nao depende de nada}
	Our choice of $r$ does not depend on $M$,$n$, $\lambda$ or $k$. This free determination of $r$ will be useful further down in the argumentation, because using the embedding $W^{1,2}(\bb{R})\hookrightarrow L^{\infty}(\bb{R})$ we will be able to obtain a uniform upper bound, in the norm of $L^{\infty}(\bb{R})$, for the sequence of solutions of the problem (\ref{problem with f_k}). Then this upper bound will naturally be transferred to also bound the sequence of solution of (\ref{main problem}).
\end{remk}
Since $\|v_M\|_{W^{1,2}}\leq r$ there is $v_0\in \bb{E}_0^1(-n,n)$ such that $v_M\rightharpoonup v_0$ in $H_0^1(-n,n)$. By the compact embedding $W^{1,2}(-n,n)\hookrightarrow L^2(-n,n)$ we conclude $v_M\rightarrow v_0$ in $L^2(-n,n)$. Our goal is to show that $v_0$ is a weak solution of (\ref{problem with f_k}). Let $\Gamma_M:V_M\rightarrow V_M^*$ and $B_M:V_M\rightarrow V_M ^*$ be defined by \begin{align}
	&\pinterno{\Gamma_M(v)}{\varphi}= \int_{-n}^{n}A(v)v'\varphi'\\
	\text{and}\;\\ \nonumber
	& \pinterno{B_M(v)}{\varphi}=\int_{-n}^{n}\left(-v+ \lambda a_1  |v|^{q-1}+|v|^{p-1}+ f_k(|v'|)+ \frac{\psi}{k}\right)\varphi.
\end{align}
Hence, $\pinterno{\Gamma_M(v_M)-B_M(v_M)}{\varphi}=0$ for all $\varphi \in V_M$.

Denoting $P_M:\bb{E}_0^1(-n,n)\rightarrow V_M$ the projection of $\bb{E}_0^1(-n,n)$ onto  $V_M$, (that is, if $u=\sum_{i=1}^{\infty}\alpha_i e_i$ then $P_M(u)=\sum_{i=1}^{M}\alpha_i e_i $) we have \begin{equation*}\pinterno{\Gamma_M(v_M)-B_M(v_M)}{v_M-P_Mv_0}=0,\end{equation*}
so
\begin{align}\label{convergence useful to prove convergence of the derivative functional form}
	&\pinterno{\Gamma_M(v_M)}{v_M-P_Mv_0}= \pinterno{B_M(v_M)}{v_M-P_Mv_0}=\\ \nonumber
	&= \int_{-n}^{n}\left(-v_M+ \lambda a_1 |v_M|^{q-1}+|v_M|^{p-1}+ f_k(|v'_M|)+ \frac{\psi}{k}\right)(v_M-P_Mv_0).
\end{align}
Letting $M\rightarrow \infty$ one can see without difficulties that $\pinterno{\Gamma_M(v_M)}{v_M-P_Mv_0}\rightarrow 0$. This convergence allows us to prove the following
\begin{lema}\label{lema-provar-convegencias}
	$v_M \rightarrow v_0$ strongly, i.e. in the norm of $H_0^1(-n,n)$.
\end{lema}
\begin{remk}
	The idea to consider the operators $\Gamma_M$ and $B_M$ was an inspiration from the arguments presented in \cite{2016_Faria}.
\end{remk}

\begin{proof}
	The limit $\|v_M-v_0\|_{L^2(-n,n)}\rightarrow 0$ has been established before, thus we will focus our efforts demonstrating the same for $\|v'_M-v'_0\|_{L^2(-n,n)}$. Let $\Phi_M,\Phi,\Psi_M,\zeta_M\in (\bb{E}_0^1(-n,n))^*$ be given by
	\begin{align}
		\Phi_M(w)= \int_{-n}^{n}A(v_M)v'_0w'\\
		\Phi(w)= \int_{-n}^{n}A(v_0)v'_0w'\\
		\Psi_M(w)= \int_{-n}^{n}A(v_M)P_Mv'_0w'\\
		\zeta_M(w)=\int_{-n}^{n}A(v_0)P_Mv'_0w'.
	\end{align}
	Then, by a straightforward calculation,  $|\Phi_M-\Phi|\rightarrow 0$,$|\Psi_M-\Phi_M|\rightarrow 0$ and $|\zeta_M-\Phi|\rightarrow 0$ in $(\bb{E}_0^1(-n,n))^*$. Thus, $|\Psi_M-\zeta_M|\rightarrow 0$ in $(\bb{E}_0^1(-n,n))^*$, since $|\Psi_M-\zeta_M|\leq |\Psi_M-\Phi_M| + |\Phi_M-\Phi| + |\Phi-\zeta_M|$. Writing $\Psi_M =(\Psi_M-\zeta_M) +\zeta_M$ yields that $\Psi_M\rightarrow \Phi$ in $(\bb{E}_0^1(-n,n))^*$. Remembering the weak convergence $v_M\rightharpoonup v_0$ one can conclude $(v_M-P_Mv_0)\rightharpoonup 0$ in $\bb{E}_0^1(-n,n)$ because for all $f\in (\bb{E}_0^1(-n,n))^*$
	\begin{align*}
		|f(v_M)-f(P_Mv_0)|\leq |f(v_M)-f(v_0)|+\|f\|\|v_0-P_Mv_0\|_{W^{1,2}}.
	\end{align*}
	Consequently, letting $M\rightarrow \infty$, $\Psi_M(v_M-P_Mv_0)\rightarrow \Phi(0)=0$. This means that
	\begin{equation}\label{convergence Luis noticed to be usefull}
		\int_{-n}^{n}A(v_M)P_Mv'_0(v'_M-P_Mv'_0)\rightarrow 0.
	\end{equation}
	Also, rewriting (\ref{convergence useful to prove convergence of the derivative functional form})
	\begin{equation}\label{convergence useful to prove convergence of the derivative integral form}
		\int_{-n}^{n}A(v_M)v'_M(v'_M-P_Mv'_0)\rightarrow 0\;\text{as}\; M\rightarrow \infty.
	\end{equation}
	Therefore, from (\ref{convergence useful to prove convergence of the derivative integral form})--(\ref{convergence Luis noticed to be usefull})
	\begin{equation}
		\int_{-n}^{n}A(v_M)(v'_M-P_Mv'_0)^2\rightarrow 0\;\text{as}\; M\rightarrow \infty.
	\end{equation}
	Since $A(x)\geq \gamma>0$ for all $x\in \bb{R}$ we conclude $\|v'_M-P_Mv'_0\|_{L^2(-n,n)}\rightarrow0$ as $M\rightarrow \infty$. Then $\|v'_M-v'_0\|_{L^2(-n,n)}\rightarrow 0$ as result of $\|v'_M-v'_0\|_{L^2(-n,n)}\leq \|v'_M-P_Mv'_0\|_{L^2(-n,n)}+\|v'_0-P_Mv'_0\|_{L^2(-n,n)}$, proving the Lemma.
\end{proof}
We know that for every $\varphi \in V_M$
\begin{align}
	&\int_{-n}^{n}A(v_M)v'_M\varphi' + \int_{-n}^{n}v_M\varphi =\\ \nonumber
	&=\int_{-n}^{n}\lambda a_1|v_M|^{q-1}\varphi + \int_{-n}^{n}|v_M|^{p-1}\varphi +\int_{-n}^{n}f_k(|v'_M|)\varphi + \int_{-n}^{n}\frac{\psi}{k}\varphi.
\end{align}
By the previous lemma, taking a subsequence, if necessary, we may assume that $v'_M(x)$ converges a.e. to $v'_0(x)$ and there exists $h\in L^2(-n,n)$ such that $|v'_M(x)|\leq h(x)$ a.e. Then notice that
\begin{align}
	\left|\int_{-n}^{n}\left(A(v_M)v'_M-A(v_0)v'_0\right)\varphi'\right| \leq \left(\int_{-n}^{n}|A(v_M)v_M'-A(v_0)v'_0|^2\right)^{1/2}\|\varphi'\|_{L^2}
\end{align}
and exists $Q>0$ such that $\|v_M\|_{\infty}< Q$ for all $M\in \mathbb{N}$, because $v_M$ converges to $v_0$ in $C^0[-n,n]$ due to the embedding $W^{1,2}(-n,n)\hookrightarrow C^0[-n,n]$. We can  suppose $Q$ big enough so that $Q>\max\{\|v_0\|_\infty+A(0),\|v_M\|_\infty+A(0)\}$ . Since
\begin{equation}
	|A(v_M(x))v_M'(x)-A(v_0(x))v'_0(x)|\rightarrow 0\;\; \text{a.e.}
\end{equation}
and
\begin{align*}
	|A(v_M(x))v_M'(x)-A(v_0(x))v'_0(x)|^2&\leq \left(|A(v_M(x))v'_M(x)|+|A(v_0(x))v'_0(x)|\right)^2 \\
	&\leq |A(v_M(x))|^2|v'_M(x)|^2 \\
	&\rquad{1}+ 2|A(v_M(x))||A(v_0(x))||v_M'(x)||v_0'(x)|\\ &\rquad{1}+|A(v_0(x))|^2|v'_0(x)|^2\\
	&\leq \tilde{A}^2Q^2h^2(x)+2\tilde{A}^2Q^2|v'_0(x)|h(x)\\
	&\rquad{1}+ \tilde{A}^2Q^2|v'_0(x)|^2
\end{align*}
almost everywhere, we conclude by (D.C.T) that
\begin{align}
	\int_{-n}^{n}A(v_M)v'_M\varphi'\rightarrow \int_{-n}^{n}A(v_0)v'_0\varphi' \quad \text{as}\; M\rightarrow \infty.
\end{align}
Also, by direct calculation, the following convergences are true

\begin{align}
	&\int_{-n}^{n}v_M\varphi \rightarrow \int_{-n}^{n}v_0\varphi\\
	&\int_{-n}^{n}\lambda a_1|v_M|^{q-1}\varphi\rightarrow \int_{-n}^{n}\lambda a_1|v_0|^{q-1}\varphi\\
	&\int_{-n}^{n}|v_M|^{p-1}\varphi\rightarrow \int_{-n}^{n}|v_0|^{p-1}\varphi \\
	&\int_{-n}^{n} f_k(|v'_M|)\varphi \rightarrow \int_{-n}^{n} f_k(|v'_0|)\varphi
\end{align}
as $M\rightarrow \infty$. Thus, for every $\varphi\in V_M$
\begin{align}
	&\int_{-n}^{n}A(v_0)v'_0\varphi' + \int_{-n}^{n}v_0\varphi =\\ \nonumber
	&=\int_{-n}^{n}\lambda a_1|v_0|^{q-1}\varphi + \int_{-n}^{n}|v_0|^{p-1}\varphi +\int_{-n}^{n}f_k(|v'_0|)\varphi + \int_{-n}^{n}\frac{\psi}{k}\varphi.
\end{align}
Furthermore, for every $u\in \bb{E}_0^1(-n,n)$, it follows that
\begin{align}
	&\int_{-n}^{n}A(v_0)v'_0P_Mu'\rightarrow \int_{-n}^{n}A(v_0)v'_0u'\\
	&\int_{-n}^{n}v_0P_Mu \rightarrow \int_{-n}^{n}v_0u\\
	&\int_{-n}^{n}\lambda a_1|v_0|^{q-1}P_Mu\rightarrow \int_{-n}^{n}\lambda a_1|v_0|^{q-1}u\\
	&\int_{-n}^{n}|v_0|^{p-1}P_Mu\rightarrow \int_{-n}^{n}|v_0|^{p-1}u\\
	&\int_{-n}^{n} f_k(|v'_0|)P_M u \rightarrow \int_{-n}^{n} f_k(|v'_0|)u
\end{align}
as $M\rightarrow \infty$. Thus, for every $u\in \bb{E}_0^1(-n,n)$
\begin{align}
	&\int_{-n}^{n}A(v_0)v'_0u' + \int_{-n}^{n}v_0u =\\ \nonumber
	&=\int_{-n}^{n}\lambda a_1|v_0|^{q-1}u + \int_{-n}^{n}|v_0|^{p-1}u +\int_{-n}^{n}f_k(|v'_0|)u + \int_{-n}^{n}\frac{\psi}{k}u.
\end{align}
So $v_0$ is an \emph{E-weak solution} of (\ref{problem with f_k}); by Lemma \ref{E-weak solution is also weak solution} $v_0$ is also a weak solution. Notice that
\begin{equation*}
	\|v_0\|_{W^{1,2}}\leq r,
\end{equation*}
and our choice of $r$ does not depend on $n,\lambda$ or $k$.
This finishes the proof of Proposition \ref{proposicao-solucao-fraca-problema-aproximado}.

In what follows we will make $k\rightarrow \infty$ thus we can consider $\psi \equiv 1$, because the term $\frac{\psi}{k}$ will converge to $0$ as $k\rightarrow \infty$.

\begin{prop}\label{proposicao-regularidade-solucao-problema-aproximado}
	The above weak solution $v_0$ satisfies:
	\begin{enumerate}
		\item There exist $\beta(L/\gamma) $ and $\hat{C}(L/\gamma,n)$, such that  $v_0\in C^{1,\beta}[-n,n]\cap C^2(-n,n)$ and \[|v_0|_{1+\beta}\leq \hat{C},\]
		where \begin{multline*}
			L>2\max\{Cr+\lambda \max\{|a_1(-n)|,|a_1(n)|\}|Cr|^{q-1}\\ +|Cr|^{p-1}+1,2C_1,A(Cr
			),\tilde{A}\};\end{multline*}
		\hyperlink{prop-candidato-nao-neg}{	\item $v_0(t)\geq 0$.}
	\end{enumerate}
\end{prop}
\begin{proof}
	To prove \textit{1}, we will use \cite[Theorem~1]{1988_Lieberman}. 
	Let $F:[-n,n]\times [-Cr,Cr]\times \mathbb{R}\rightarrow \bb{R}$ be defined by $F(x,z,p)= A(z)p$, where $C$ is the embedding constant for $W^{1,2}(\mathbb{R})\hookrightarrow L^{\infty}(\bb{R}),$ and $B(x,z,p)=z-(\lambda a_1(x)|z|^{q-1}+|z|^{p-1}+f_k(|p|)+\frac{1}{k})$ be defined in the same domain. Then, problem (\ref{problem with f_k}) may be rewritten as

	\[div_xF(x,u(x),u'(x))+ B(x,u(x),u'(x))=0.\]

	In order to use \cite[Theorem~1]{1988_Lieberman} 
	we must verify the existence of nonnegative constants $l,L,M_0,m, \kappa$ with $l\leq L$ such that
	\begin{align}
		&\frac{\partial F}{\partial p}(x,z,p) \xi^2 \geq l(\kappa+|p|)^m\xi^2,\label{cota-a-liberman}\\
		&\left|\frac{\partial F}{\partial p}(x,z,p)\right|\leq L(\kappa+|p|)^m,\label{cota-b-liberman}\\
		&|F(x,z,p)-F(y,w,p)|\leq L(1+|p|)^{m+1}\cdot|z-w|,\label{cota-c-liberman}\\
		&|B(x,z,p)|\leq L(1+|p|)^{m+2},\label{cota-d-liberman}
	\end{align}
	for all $(x,z,p)\in \{-n,n\}\times [-M_0,M_0]\times \bb{R}$, $w\in [-M_0,M_0]$ and $\xi \in \bb{R}$. Since $\frac{\partial F}{\partial p}(x,z,p)= A(z)$, inequality (\ref{cota-a-liberman}) follows from $A(z)\xi^2\geq \gamma \xi^2$, that is, $l=\gamma$.

	To prove the remaining inequalities take $M_0=Cr,$ \[T> \max\{Cr+\lambda \max\{|a_1(-n)|,|a_1(n)|\}|Cr|^{q-1}+|Cr|^{p-1}+1,2C_1,A(Cr
	),\tilde{A}\},\]$L=2T,\kappa =0$ and $m=0$, where $\tilde{A}$ is the Lipchitz constant of $A$. Then :
	\begin{itemize}

		\item[(\ref{cota-b-liberman}) ]\[\left|\frac{\partial F}{\partial p}(x,z,p)\right|= A(z)\leq A(Cr)<L;\]
		\item[(\ref{cota-c-liberman})]\[|F(x,z,p)-F(y,w,p)|=|A(z)p-A(w)p|\leq \tilde{A}|p||z-w|\leq L(1+|p|)|z-w|;\]
		\item[(\ref{cota-d-liberman})] \begin{align}\label{estimativa-B-liberman}
			|B(x,z,p)|&=|z-(\lambda a_1(x)|z|^{q-1}+|z|^{p-1}+f_k(|p|)+\frac{1}{k})|\\\nonumber
			&\leq Cr + \lambda \max\{|a_1(-n)|,|a_1(n)|\}|Cr|^{q-1} \\\nonumber
			&\;\;\rquad{2}+ |Cr|^{p-1}+1/k + C_1(1+|p|^{\theta-1})\\\nonumber
			&\leq T+C_1(1+(1+|p|)^{\theta-1})\\\nonumber
			&\leq T +2 C_1(1+|p|)^2\\\nonumber
			&\leq T(1+(1+|p|)^2)\\\nonumber
			&\leq 2T(1+|p|)^2 = L(1+|p|)^2.\\\nonumber
		\end{align}
	\end{itemize}
	Therefore, by \cite[Theorem~1]{1988_Lieberman}
	there exists $\beta \in (0,1)$ and a constant $\hat{C}$, independent of $k$, such that $v_0\in C^{1,\beta}([-n,n])$ and
	\begin{equation}
		|v_0|_{1+\beta}\leq \hat{C}.\label{estimativa-liberman-norma}\end{equation}
	It follows from \cite[pag.~317, Chap.~6,  Theorem~4]{2010_Evans_BOOK} 
	that $v_0 \in W^{2,2}(-n,n)$ and since $v_0$ is a weak solution of (\ref{problem with f_k}) we have
	\begin{equation}\label{cara-da-segunda-derivada}
		v_0''=\frac{v_0-\lambda a_1|v_0|^{q-1}-|v_0|^{p-1}-f_k(|v_0'|)-1/k - A'(v_0)|v_0'|^2}{A(v_0)}
	\end{equation}
	showing that $v_0''$ is continuous, thus $v_0\in C^2(-n,n).$

	\hypertarget{prop-candidato-nao-neg}{	To prove that $v_0(t)\geq 0$ for all $t\in (-n,n)$} we first notice that $v_0^{-}(t)=\max\{0,-v_0(t)\}\in H_0^1(-n,n)$.  Using $v_0^-(t)$ as a test function in the definition of weak solution provides
	\begin{multline}
		-\int_{-n}^{n}A(v_0)|v_0^{-}|^2 - \int_{-n}^{n}|v_0^{-}|^2 = \int_{-n}^{n}\lambda a_1|v_0|^{q-1}v_0^{-}+\int_{-n}^{n}|v_0|^{p-1}v_0^{-}\\+\int_{-n}^{n}f_k(|v_0'|)v_0^{-} + \int_{-n}^{n}\frac{1}{k} v_0^{-}.
	\end{multline}
	Then $-\gamma\|v_0^{-}\|_{W^{1,2}}^2\geq 0$, thus $\|v_0^{-}\|_{W^{1,2}}=0$ implying $v_0^{-}\equiv 0$ a.e. Since $v_0$ is continuous, $v_0(t)\geq 0$ for all $t\in (-n,n).$ This finishes the proof of Proposition \ref{proposicao-regularidade-solucao-problema-aproximado}.

\end{proof}

Thus, by Proposition \ref{proposicao-solucao-fraca-problema-aproximado} and Proposition \ref{proposicao-regularidade-solucao-problema-aproximado} we obtain the proof of Theorem \ref{teo-da-solucao-do-problema-aproximado}.

\subsection{Constructing a solution to problem (\ref{main problem})}
\label{subsec: constructing a solution to problem Pn}

Let $v_k$ be the (strong) solution of problem (\ref{problem with f_k}) -- obtained just above -- with $k$ varying. By the previous constructions, we have that $\normaSobolev{v_k}{-n,n} \leq r$ independent of $k$, as noticed in Remark \ref{observacao r nao depende de nada}. Then there exists $u_n\in H_0^1(-n,n)$, $\normaSobolev{u_n}{-n,n}\leq r$, so that $v_k$ has a subsequence  converging weakly in $H^1_0(-n,n)$ to $u_n$. From now on $v_k$ will denote this subsequence. Since the function
\begin{align*}
	H_0^1(-n,n)\ni w \mapsto \int_{-n}^{n} A(u_n)u'_nw'
\end{align*}
belongs to $(H_0^1(-n,n))^*$, we have -- by the weak convergence -- that
\[\int_{-n}^{n} A(u_n)u'_n(v_k-u_n)'\rightarrow 0\; \mbox{ as } k\rightarrow \infty.\]
This convergence will be useful in our next task: to prove that $v_k\rightarrow u_n$ strongly in $H^1_0(-n,n)$.
\begin{lema}
	The following convergence is true
	\[\int_{-n}^{n}A(u_n)v'_k(v_k-u_n)'\rightarrow 0\mbox{ as }k\rightarrow \infty.\]
\end{lema}
\begin{proof}
	We might write
	\begin{align*}
		\int_{-n}^{n}A(u_n)v'_k(v_k-u_n)' &= \int_{-n}^{n}[A(u_n)-A(v_k)+A(v_k)]v'_k(v'_k-u'_n)\\
		&\negarquad{2}= \underbrace{\int_{-n}^{n}[A(u_n)-A(v_k)]v'_k(v'_k-u'_n)}_{I_1} +\underbrace{ \int_{-n}^{n}A(v_k)v'_k(v'_k-u'_n)}_{I_2}
	\end{align*}
	and analyze $I_1$ and $I_2$ separately.

	\paragraph{Analysis of $I_2$ } By the weak formulation of (\ref{problem with f_k})
	\begin{align*}
		\int_{-n}^{n}A(v_k)v'_k(v'_k-u'_n)=\negarquad{10} &\nonumber\\
		&\underbrace{\int_{-n}^{n}-v_k(v_k-u_n)	+ \lambda a_1|v_k|^{q-1}(v_k-u_n) +|v_k|^{p-1}(v_k-u_n)+ \frac{(v_k-u_n)}{k}}_{E_1}
		\\
		&\underbrace{+ \int_{-n}^{n} f_k(|v'_k|)(v_k-u_n)}_{E_2} 
		.
	\end{align*}
	Since we have compact injection of $H_0^1(-n,n)$ onto $L^2(-n,n)$, the weak convergence of $v_k$ to $u_n$ in $H^1_0(-n,n)$ implies $\|v_k-u_n\|_{L^2}\rightarrow 0$. Thus, it is straightforward to see that ($E_1$) 
	converges to $0$ as $k\rightarrow \infty$. Remains to verify what happens with ($E_2$) 
	in the limit. We have that
	\begin{align*}
		\int_{-n}^{n}f_k(|v_k'|)(v_k-u_n)\leq \int_{-n}^{n}C_1(|v_k'|^{\theta-1}+|v'_k|)|v_k-u_n|
	\end{align*}
	by Lemma \ref{upper-estimative for the seq f_k}. Using Proposition \ref{proposicao-regularidade-solucao-problema-aproximado}, item \emph{1}, that is, the estimation $|v_k|_{1,\beta}\leq \hat{C}$ which is independent of $k$, we have that
	\[|v_k'|^{\theta-1}+|v_k'|\leq (\hat{C})^{\theta-1}+\hat{C}.\]
	Then,
	\begin{align*}
		\int_{-n}^{n}f_k(|v_k'|)(v_k-u_n)&\leq C_1[(\hat{C})^{\theta-1}+\hat{C}]\int_{-n}^{n}|v_k-u_n|\\
		&\leq \underbrace{(2n)^{1/2}C_1[(\hat{C})^{\theta-1}+\hat{C}]\|v_k-u_n\|_{L^2}}_{\rightarrow 0 \mbox{ as }k\rightarrow \infty}.
	\end{align*}
	Thus, $\lim_{k\rightarrow \infty} I_2(k)=0$.

	\paragraph{Analysis of $I_1$ } We also have that $\lim_{k\rightarrow \infty} I_1(k)=0$, as one can see through
	\begin{align*}
		\left|\int_{-n}^{n}[A(u_n)-A(v_k)]v'_k(v'_k-u'_n)\right| &\leq \int_{-n}^{n}|A(u_n)-A(v_k)||v'_k||v'_k-u'_n|\\
		&\leq \hat{C} \int_{-n}^{n}|A(u_n)-A(v_k)||v'_k-u'_n|\\
		&\leq\hat{C}\tilde{A}\int_{-n}^{n}|u_n-v_k||v'_k-u'_n|\\
		&\leq \hat{C}\tilde{A} \|u_n-v_k\|_{L^2}\|v'_k-u'_n\|_{L^2}\\
		&\leq \hat{C}\tilde{A}2r\|u_n-v_k\|_{L^2}.
	\end{align*}
	Proving the Lemma.$\;\;\;\rquad{30}$\qedsymbol

	Thus,

	\begin{align}
		&\int_{-n}^{n} A(u_n)u'_n(v_k-u_n)'\rightarrow 0\; \mbox{ as } k\rightarrow \infty\label{boa convergencia1}\\
		&\int_{-n}^{n}A(u_n)v'_k(v_k-u_n)'\rightarrow 0\mbox{ as }k\rightarrow \infty\label{boa convergencia2}.
	\end{align}
	Subtracting (\ref{boa convergencia2}) from (\ref{boa convergencia1}) we have
	\begin{equation}
		\int_{-n}^{n}A(u_n)(v'_k-u'_n)^2 \rightarrow 0 \mbox{ as } k\rightarrow \infty
	\end{equation}
	implying that $v'_k\rightarrow u'_n$ in $L^2(-n,n)$, since $\gamma$ is a uniform lower-bound for $A$. Hence $v_k\rightarrow u_n$ in $H^1_0(-n,n)$.
\end{proof}
\begin{remk}
	Since $v_k\rightarrow u_n$ in $H^1_0(-n,n)$ we conclude that $u_n$ is also an even function; due to the embedding $W^{1,2}(-n,n)\hookrightarrow C[-n,n]$.
\end{remk}

\begin{prop}\label{proposicao-propiedades-un}
	The function $u_n$ satisfies:
	\begin{enumerate}
		\item $u_n$ is strictly positive in $(-n,n)$;
		\hyperlink{Item2-prop-eh-solucao-no-intervalo}{\item $u_n$ is a solution to} (\ref{main problem}).
	\end{enumerate}
\end{prop}

\begin{proof}
	$\quad$

	\paragraph{Item 1. }Let $\tilde{a}:= \inf_{x\in [-n,n]} a_1(x)$. We will divide our argument into two cases:
	\begin{remk}
		This division of cases is a geometric argument that we borrowed from \cite{2011_Alves}.
	\end{remk}
	\begin{center}
		\textbf{Case 1.} There exists a subsequence $(v_{k_i})_{i\in \bb{N}}$ of $(v_k)$ such that $v'_{k_i}\geq 0$ in $(-n,0)$ for all $i$.
	\end{center}

	Consider the problem
	\begin{equation}
		\begin{cases}
			-(A(u)u')'+u = \lambda \tilde{a}|u|^{q-1}\quad \mbox{in }(-n,n)\\
			u>0 \;\;\;\rquad{9} \mbox{ in }(-n,n)\\
			u(-n)=u(n)=0.
		\end{cases}\label{equacao autofuncao por baixo}
	\end{equation}
	Since $v_{k_i}'\geq 0$ in $(-n,0)$ we get that $v_{k_i}>0$ in $(-n,0)$, because -- due to Proposition \ref{proposicao-regularidade-solucao-problema-aproximado} -- $v_{k_i}$ is an even solution of (\ref{problem with f_k}), thus \emph{it can not be identically zero in an interval} and  $v_{k_i}\geq0$ ; i.e., supposing the existence of $x_i\in (-n,0)$ such that $v_{k_i}(x_i)=0$ implies the existence of $y_i\in (-n,0)$ such that $v_{k_i}'(y_i)<0$, which would be a contradiction. Thus, we see that $v_{k_i}$ is a sup-solution for this equation. Let $\phi_1$ be an even and positive eigenfunction  for the eigenvalue problem
	\begin{equation}
		\begin{cases}
			-u''=\lambda_1 u \quad \mbox{in } (-n,n)\\
			u(-n)=u(n)=0
		\end{cases}
	\end{equation}
	where $\lambda_1= \frac{\pi^2}{(2n)^2}$. Thus, choosing $\tau$ such that
	\[\frac{\tau^{2-q}(1+\gamma\lambda_1)}{\lambda\tilde{a}}\leq \phi_1^{q-2}\]
	we have that $\tau\phi_1$ is as sub-solution of (\ref{equacao autofuncao por baixo}). By Theorem \ref{criterio-de-comparação} 
	\[v_{k_i}(t)\geq \tau \phi_1(t)\quad \forall t\in (-n,n),\]
	therefore, in the limit,
	\[u_n(t)\geq \tau\phi_1(t)>0\quad \forall t\in (-n,n).\]

	\begin{center}
		\textbf{Case 2.}
		There exists a subsequence $(v_{k_i})_{i\in \bb{N}}$ of $(v_k)$ and there exists a sequence $(z_i)_{i\in \mathbb{N}}\subset(-n,0)$ such that $v'_{k_i}(z_i)<0$.
	\end{center}

	\begin{remk}
		Although the geometric argument is an inspiration from \cite{2011_Alves}, we still need to adjust it to our necessity. Lemma \ref{lema-caso2-limita-por-baixo} is one such adjustment.
	\end{remk}

	\begin{lema}
		\label{lema-caso2-limita-por-baixo}
		Let $x\in (-n,n)$ such that $v''_{k_i}(x)\geq 0$, then $v_{k_i}(x)> (\lambda\tilde{a})^{\frac{1}{2-q}}$.
	\end{lema}
	\begin{proof}
		Since $v_{k_i}$ is a solution for the problem (\ref{problem with f_k}), with $k=k_i$, for all $t$
		\begin{align*}
			-A'(v_{k_i}(t))|v'_{k_i}(t)|^2-A(v_{k_i}(t))v''_{k_i}(t)+v_{k_i}(t)& \\
			&\hspace{-93pt}=\lambda a_1(t)|v_{k_i}(t)|^{q-1} + |v_{k_i}(t)|^{p-1}+ f_{k_i}(|v_{k_i}'(t)|)+ \frac{1}{k_i}\\
			&\hspace{-92pt}\geq \lambda a_1(t)|v_{k_i}(t)|^{q-1} + |v_{k_i}(t)|^{p-1} + \frac{1}{k_i}.
		\end{align*}
	Here we used that $sign(f_{k_i}(s))= sign(s)$, thus $f_{k_i}(|v'_{k_i}|)\geq 0$. Using that $|v_{k_i}|^{p-1}\geq 0$ and $a_1(t)\geq \tilde{a}>0$ we obtain:
		\begin{equation}
			-A'(v_{k_i}(t))|v'_{k_i}(t)|^2-A(v_{k_i}(t))v''_{k_i}(t)+v_{k_i}(t)\geq \lambda\tilde{a}|v_{k_i}(t)|^{q-1}+\frac{1}{k_i}.
		\end{equation}
		Then, with $t=x$,
		\begin{align*}
			-A(v_{k_i}(x))v''_{k_i}(x)&\geq \lambda\tilde{a}|v_{k_i}(x)|^{q-1}	-v_{k_i}(x)+A'(v_{k_i}(x))|v'_{k_i}(x)|^2+ \frac{1}{k_i}\\
			&> \lambda\tilde{a}|v_{k_i}(x)|^{q-1}	-v_{k_i}(x).
		\end{align*}
	Where, in the last inequality, we used that $A$ is non-decreasing, $|v'_{k_i}|\geq 0$ and $1/k_i>0$. Notice that the resulting estimation is strict because $1/k_i>0$. By hypotheses $v''_{k_i}(x)\geq 0$, then $-A(v_{k_i}(x))v''_{k_i}(x)\leq 0$. Using the previous inequality,
		\[v_{k_i}(x)>\lambda\tilde{a}|v_{k_i}(x)|^{q-1},\]
		thus $v_{k_i}(x)\neq 0$ and  $v_{k_i}(x)>(\lambda\tilde{a})^{\frac{1}{2-q}}$.
	\end{proof}
	Now, in order to use this lemma, we ought to find a $x_i\in(-n,n)$ such that $v''_{k_i}(x_i)\geq 0$.	Using the fact that $v_{k_i}$ is even and $v'_{k_i}(z_i)<0$ we have that $v'_{k_i}(-z_i)>0$. Let $x_i=\min_{x\in[z_i,-z_i]}v_{k_i}(x)$ and notice that $x_i\neq z_i$ and $x_i\neq -z_i$; indeed, there exist $\delta>0$ such that, if $x\in (z_i,z_i+\delta)\cup (-z_i-\delta,-z_i)$, then $v_{k_i}(x)< v_{k_i}(z_i)=v_{k_i}(-z_i).$ Hence $x_i\in (z_i,-z_i)$ and $v'_{k_i}(x_i)=0$; therefore $v''_{k_i}(x_i)$ must be \emph{greater or equal than 0}, because if $v''_{k_i}(x_i)<0$ there would be $\xi>0$ such that, for $x\in(x_i,x_i + \xi)\subset (z_i,-z_i)$, $v'_{k_i}(x)<0$; and for this neighborhood $v_{k_i}(x)<v_{k_i}(x_i)$ -- a contradiction with the minimality of $x_i$. Thus,  $v''_{k_i}(x_i)\geq 0$. By Lemma \ref{lema-caso2-limita-por-baixo}, we obtain 
	
	for all $i$
	\[v_{k_i}(x_i)> (\lambda\tilde{a})^{\frac{1}{2-q}}.\]
	From the compacity of $[-n,n]$, there exist $x_0\in [-n,n]$ such that $x_i\rightarrow x_0$ when $i\rightarrow \infty$; taking a subsequence if necessary. Then
	\[u_n(x_0)=\lim_{i\rightarrow \infty}v_{k_i}(x_i)\geq (\lambda\tilde{a})^{\frac{1}{2-q}}>0.\]
	Finally, we will conclude \emph{item 1} showing that, also in this case, $u_n$ is strictly positive in $(-n,n)$.

	Suppose by contradiction that there exists $y\in (-n,n)$ such that $u_n(y)=0$. Let $(d,s)\subset (-n,n)$ be the biggest interval containing $y$ satisfying the property: if $x\in (d,s)$ then $u_n(x)< \frac{(\lambda\tilde{a})^{\frac{1}{2-q}}}{2}$. Since $u_n(x_0)=u_n(-x_0)> \frac{(\lambda \tilde{a})^{\frac{1}{2-q}}  }{2}$ we have that $d\neq -n$ or $s\neq n$. Thus we can suppose without loss of generality that $d>-n$, because on the contrary, we would apply our following arguments using the interval $(d',s')$, where $d'=-s$ and $s'=-d$, and the point $y'=-y$.

	By the maximality  of $(d,s)$ and the continuity of $u_n$ we have that $u_n(d)=\frac{(\lambda\tilde{a})^{\frac{1}{2-q}}}{2}$. Since $u_n(x)<\frac{(\lambda\tilde{a})^{\frac{1}{2-q}}}{2}$ for all $x\in (d,s)$ and $v_{k_i}$ converges uniformly to $u_n$, there exist $i_1\in \bb{N}$ such that, for $i>i_1$ and $x\in (d,s)$,
	\[v_{k_i}(x)<(\lambda\tilde{a})^{\frac{1}{2-q}}.\]
	Then, by Lemma \ref{lema-caso2-limita-por-baixo} $v''_{k_i}(x)<0$. Using that $u_n(d)=\frac{(\lambda\tilde{a})^{\frac{1}{2-q}}}{2}$, there exist $i_2\in \bb{N}$ such that $i>i_2$ implies
	\[v_{k_i}(d)>\frac{(\lambda\tilde{a})^{\frac{1}{2-q}}}{4}.\]
	Let $i_0>\max\{i_1,i_2\}$ and define $f:(d,s)\rightarrow \bb{R}$ by
	\[f(x)=\frac{(\lambda\tilde{a})^{\frac{1}{2-q}}}{4}\cdot\frac{x-s}{d-s}.\]
	We have that $f(d)=\frac{(\lambda\tilde{a})^{\frac{1}{2-q}}}{4}$ and $f(s)=0$. Let $U_i(x)=v_{k_i}(x)-f(x)$ for $i\geq i_0$, then
	\begin{equation}
		\begin{cases}
			U''_i(x)<0,\quad \mbox{for }x\in (d,s)\\
			U_i(d)>0,U_i(s)=v_{k_i}(s)\geq0.
		\end{cases}
	\end{equation}
	By the maximum principle, the minimum of $U_i$ is reached on the border of the interval $(d,s)$, implying that $U_i(x)>0$ for all $x\in (d,s)$, that is, $v_{k_i}(x)>f(x)$ for all $x\in (d,s)$ and $i\geq i_0$. Thus, taking $x=y$ and making $i\rightarrow \infty$, we obtain
	\[u_n(y)\geq f(y)>0,\]
	which is a contradiction.

	\paragraph{Item 2. } \hypertarget{Item2-prop-eh-solucao-no-intervalo}{Since the estimation} from Proposition \ref{proposicao-regularidade-solucao-problema-aproximado} item \emph{1} 
	holds, that is,
	\[|v_k|_{1+\beta}\leq \hat{C}\]
	for all $k\in \mathbb{N}$; and, for all $1<\alpha<\beta$, we have compact embedding $C^{1,\beta}[-n,n]\hookrightarrow C^{1,\alpha}[-n,n]$, we may assume -- taking a subsequence, if necessary -- that there exist $\widetilde{u_n}\in C^{1,\alpha}[-n,n]$ such that $v_k\rightarrow \widetilde{u_n}$ in $C^{1,\alpha}[-n,n]$ as $k\rightarrow \infty$. Thus,
	\begin{align*}
		&v_k\rightarrow u_n \mbox{ in } C^0[-n,n] \mbox{ as } k\rightarrow \infty\\
		&v_k\rightarrow \widetilde{u_n} \mbox{ in } C^{1,\alpha}[-n,n] \mbox{ as } k\rightarrow \infty.
	\end{align*}
	Then for all $x\in [-n,n]$ we have
	\[u_n(x)=\lim_{k\rightarrow \infty}v_k(x)=\widetilde{u_n}(x),\]
	i.e., $u_n=\widetilde{u_n}\in C^{1,\alpha}[-n,n]$.

	Considering the definition of weak solution, for all $\varphi \in H^1_0(-n,n)$
	\begin{multline*}
		\int_{-n}^{n}A(v_k)v'_k\varphi' + \int_{-n}^{n}v_k\varphi = \int_{-n}^{n}(\lambda a_1 |v_k|^{q-1}+|v_k|^{p-1})\varphi  \\ + \int_{-n}^{n} f_k(|v'_k|)\varphi + \int_{-n}^{n} \frac{\varphi}{k}.
	\end{multline*}
	By (D.C.T), it is straightforward to see that the following convergences are true:
	\begin{align*}
		&\int_{-n}^{n}A(v_k)v'_k\varphi' \rightarrow 	\int_{-n}^{n}A(u_n)u'_n\varphi',\\
		&\int_{-n}^{n}v_k\varphi\rightarrow
		\int_{-n}^{n}u_n\varphi,\\ &\int_{-n}^{n}(\lambda a_1 |v_k|^{q-1}+|v_k|^{p-1})\varphi\rightarrow \int_{-n}^{n}(\lambda a_1 |u_n|^{q-1}+|u_n|^{p-1})\varphi,\\
		&\int_{-n}^{n} \frac{\varphi}{k} \rightarrow 0,
	\end{align*}
	as $k\rightarrow \infty$. Let us examine the remaining integral. First notice that $f_k(|v'_k|)$ converges uniformly to $g(|u'_n|)$; indeed, given $\epsilon>0$ there exists $k_0\in \mathbb{N}$ such that $k>k_0$ implies
	\begin{equation}
		||v'_k(x)|-|u'_n(x)||< \epsilon \quad\forall x\in [-n,n].
	\end{equation}
	Also there exist $0<\delta_k<\epsilon$ such that if $|x-y|<\delta_k$, then
	\begin{equation}
		|f_k(x)-f_k(y)|<\frac{\epsilon}{2};
	\end{equation}
	thus for $k>k_0$
	\begin{equation}
		|f_k(|v'_k(x)|)-f_k(|u'_n(x)|)|<\frac{\epsilon}{2} \quad \forall x \in [-n,n].
	\end{equation}
	In the perspective of Theorem \ref{teo of sequence f_k}, $f_k$ converges to $g$ uniformly in bounded sets; since $\|u'_n\|_{\infty}\leq \hat{C}$, for $x\in [-\hat{C},\hat{C}]$ there exist $k_1\in \mathbb{N}$ such that $k>k_1$ implies
	\begin{equation}
		|f_k(x)-g(x)|<\frac{\epsilon}{2}\quad \forall x\in [-\hat{C},\hat{C}]
	\end{equation}
	and with all these ingredients we obtain the uniform convergence, because for $k>\max\{k_0,k_1\}$
	\begin{align*}
		|f_k(|v'_k(x)|)-g(|u'_n(x)|)|&\leq |f_k(|v'_k(x)|)-f_k(|u'_n(x)|)|\\
		& \rquad{2}+|f_k(|u'_n(x)|)-g(|u'_n(x)|)|\\
		&<\epsilon \quad \forall x\in [-n,n].
	\end{align*}
	Thus, by (D.C.T)
	\[\int_{-n}^{n}f_k(|v'_k|)\varphi\rightarrow \int_{-n}^{n}g(|u'_n|)\varphi\]
	as $k\rightarrow \infty$. All these convergences together show that $u_n$ is a weak solution for the problem (\ref{main problem}).

	From \cite[Pag.~317, Chap.~6, Theorem~4]{2010_Evans_BOOK} 
	we conclude that $u_n\in W^{2,2}(-n,n)$; and similarly, to the argument showed in equation (\ref{cara-da-segunda-derivada}) we obtain that $u''_n\in C^0(-n,n)$. Thus, $u_n$ is a strong solution to the problem (\ref{main problem}).

\end{proof}

\section{Asymptotic solution to problem $(P_n)$}
\label{sec: asympotic solution to problem Pn}
In this section, we will briefly study the solution's behavior of problem (\ref{main problem}) when $\lambda\rightarrow 0$ or $\lambda\rightarrow \lambda^*$. Our main result is the following:
\begin{teo}
	Denoting by $u_\lambda$ the solution of (\ref{main problem}):
	\begin{enumerate}
		\item As $\lambda\rightarrow 0$, we get that $\|u_\lambda\|_{W^{1,2}}\rightarrow 0$;
		\item one can take $\lambda=\lambda^*$ and still obtain a solution to problem (\ref{main problem}).
	\end{enumerate}
\end{teo}
\begin{proof}

	Let $v_k$ be the strong solution of problem (\ref{problem with f_k}), with $\psi\equiv1$. In the previous section, we proved -- among other things -- that $v_k\rightarrow u_\lambda$ as $k\rightarrow \infty$, (taking a subsequence, if necessary). Using $v_k$ as test function in the Definition \ref{definicao-sol-fraca-problema-com-fk}, we get:

	\begin{align*}
		\int_{-n}^{n}A(v_k)|v'_k|^2 + \int_{-n}^{n}|v_k|^2 = \int_{-n}^{n}\lambda a_1 |v_k|^{q-1}v_k + \int_{-n}^{n}|v_k|^{p-1}v_k\\
		+\int_{-n}^{n}f_k(|v'_k|)v_k + \int_{-n}^{n}\frac{v_k}{k}.
	\end{align*}
	Thus, following the estimations done in Proposition \ref{proposicao-solucao-fraca-problema-aproximado}, we can estimate these integrals to obtain
	\begin{align*}
		\gamma \|v_k\|^2_{W^{1,2}}\leq \lambda C_2 \|v_k\|^q_{W^{1,2}} + C^{p-2}\|v_k\|^p_{W^{1,2}} + C_1\max\{C^{\theta-2},C\}\|v_k\|^\theta_{W^{1,2}}\\
		+\left(\frac{C_1(2n)^{1/2}}{k}+ \frac{(2n)^{1/2}}{k}\right)\|v_k\|_{W^{1,2}}.
	\end{align*}
	Rearranging we obtain
	\begin{multline}
		\|v_k\|^2_{W^{1,2}}\left(\gamma-C^{p-2}\|v_k\|^{p-2}_{W^{1,2}}-C_1\max\{C^{\theta-2},C\}\|v_k\|^{\theta-2}_{W^{1,2}}\right)\leq \\ \lambda C_2 \|v_k\|^q_{W^{1,2}} + 	\left(\frac{C_1(2n)^{1/2}}{k}+ \frac{(2n)^{1/2}}{k}\right)\|v_k\|_{W^{1,2}}.
	\end{multline}
	Notice that $\|v_k\|_{W^{1,2}}\leq r$ independent of $k$ and

	\[r< \min\left\{\left(\frac{\gamma}{4C^{p-2}}\right)^{1/(p-2)},\left(\frac{\gamma}{4C_1\max\{C^{\theta-2},C\}}\right)^{1/(\theta-2)}\right\}.\]

	Thus,
	\begin{align}
		\|v_k\|^2_{W^{1,2}} \leq \frac{2}{\gamma}\left(\lambda C_2r^q + \left(\frac{C_1(2n)^{1/2}}{k}+ \frac{(2n)^{1/2}}{k}\right) r \right).
	\end{align}
	Making $k\rightarrow \infty$ we end up with
	\begin{equation}\label{limitacao-por-cima-usando-lambda}
		\|u_\lambda\|^2_{W^{1,2}}\leq \lambda \cdot \left(\frac{2C_2r^q}{\gamma}\right).
	\end{equation}
	Then, as $\lambda\rightarrow 0$ we see that $	\|u_\lambda\|^2_{W^{1,2}}\rightarrow 0$. This proves item \textit{1}.

	To prove item \textit{2}, one can take a sequence $(\lambda_n)$ in $(0,\lambda^*)$ such that $\lambda_n\rightarrow \lambda^*$. Noticing that $\|u_{\lambda_n}\|_{W^{1,2}}\leq r$ independent of $\lambda_n$, one can obtain a candidate $u_{\lambda^*}\in H^1_0(-n,n)$ such that $u_{\lambda_n}$ converges weakly to $u_{\lambda^*}$ in $ H^1_0(-n,n)$. Then, following similar argumentation as exposed in Section \ref{sec:Solution in a bounded interval}, one can prove that $u_{\lambda^*}$ is an even, positive solution to (\ref{main problem}) with $\lambda= \lambda^*$.
\end{proof}
\begin{remk}
	Let
	\[\tilde{r} = \min\left\{ r, \lambda^{1/2}\left(\frac{2C_2 r^q}{\gamma}\right)^{1/2}   \right\}.\]
	Then, by inequality (\ref{limitacao-por-cima-usando-lambda}) and the estimations from Proposition \ref{proposicao-solucao-fraca-problema-aproximado}, we get that $\|u_\lambda\|_{W^{1,2}}\leq \tilde{r}$. Notice that $\tilde{r}\rightarrow 0$ as $\lambda\rightarrow0$.
\end{remk}

\begin{remk}
	For now, on we will resume the previous notation, that is, we will call the strong solution of problem (\ref{main problem}) by $u_n$. This will be useful in the next section.
\end{remk}

\section{Solution in $\bb{R}$}
\label{secao: solucao na reta toda}
To obtain the homoclinic solution, we can proceed similarly as in \cite{2011_Alves}. However, we present a slightly different approach.

To obtain a solution defined in $\bb{R}$ we will utilize a subsequence construction wrapping it up with the arguments presented in Section \ref{sec:Solution in a bounded interval}. The reader should notice that the notation ``$u_n$'' used for the solution of (\ref{main problem}) -- previously obtained -- in $(-n,n)$ is not accidental: extending $u_n$ by zero out of $[-n,n]$ we obtain a sequence $(u_n)$ in $H^1(\bb{R})$. Throughout this section, we will use $u_n$ to denote the solution ``$u_n$'' and its extension. Also, one can see that $\|u_n\|_{H^1(\bb{R})}=\normaSobolev{u_n}{-n,n}\leq \tilde{r}$ for all $n$.

Let $K_1=[-1,1]$; then for all $n\geq 1$ we have that $u^1_n:=u_n|_{K_1}$ is well defined and $u^1_n\in H^1(-1,1)$. By the limitation $\normaSobolev{u^1_n}{-1,1}\leq \tilde{r}$ there exists a subsequence $u_{n,1}$ and $s_1\in H^1(-1,1)$  such that $u_{n,1}\rightharpoonup s_1$ in $H^1(-1,1)$. Notice that the compact injection $H^1(-1,1)\hookrightarrow C^0[-1,1]$ implies that, passing to a subsequence,  $u_{n,1}\rightarrow s_1$ in $C^0[-1,1]$.

Let $K_2=[-2,2]$. Taking $n$ in the set of indices of the subsequence $u_{n,1}$, for $n\geq 2$ we have that $u^2_n:=u_n|_{K_2}$ is well defined and $\normaSobolev{u^2_n}{-2,2}\leq \tilde{r}$. Thus, there exists a subsequence $u_{n,2}$ of $u^2_n$ and $s_2\in H^1(-2,2)$ such that $u_{n,2}\rightharpoonup s_2$ in $H^1(-2,2)$.

Repeating the same argument, by induction we get that for all $j\in \bb{N}$ there exists a subsequence $u_{n,j}$ of $u_{n,j-1}$ and $s_j\in H^1(-j,j)$ such that $u_{n,j}\rightharpoonup s_j$ in $H^1(-j,j)$. Notice that $\|s_j\|_{H^1(-j,j)}\leq \tilde{r}$ for all $j\in \mathbb{N}$.

\begin{remk}
	$u_{n,j}$ is the subsequence of $u_n$ that converges weakly in $H^1(-j,j)$ to $s_j$. As mentioned, this weak convergence implies convergence in $C^0[-j,j]$ which gives us, in particular, that $s_j$ is an even function, since $u_{n,j}$ is even for all $n\in \mathbb{N}$.
\end{remk}

\begin{lema}\label{boa def das extensoes na reta}
	$s_j|_{[1-j,j-1]}= s_{j-1}$
\end{lema}
\begin{proof}
	Given $x\in [1-j,j-1]$ we have that
	\[s_{j-1}(x)=\lim_{n\rightarrow\infty}u_{n,j}(x)= s_j(x),\]
	because $u_{n,j}$ is a subsequence of $u_{n,j-1}$.
\end{proof}
Fix $j\in \mathbb{N}$; from here and forward we will focus our attention on proving that, in fact, $s_j$ \emph{is smooth and positive.}

 Define $W:[-n,n]\times [-Cr,Cr]\times \bb{R}\longrightarrow \bb{R}$ by
\[W(x,z,p)= z-(\lambda a_1(x)|z|^{q-1}+|z|^{p-1}+ g(|p|))\]
one can see that the estimation (\ref{estimativa-B-liberman}) also holds, with $W$ taking part as $B$ (remember that $C$ is the constant for the embedding $W^{1,2}(\bb{R})\hookrightarrow L^{\infty}(\bb{R}))$. Then, for any $n\geq j$, by Theorem 1 from \cite{1988_Lieberman} there exist $\hat{C}(j)>0$ and $0<\beta\leq 1$ such that
\[|u_n|_{1+\beta}\leq \hat{C}(j)\mbox{ in }C^{1,\beta}[-j,j].\]
Taking $0<\alpha<\beta\leq 1$ we get (see the argumentation on \hyperlink{Item2-prop-eh-solucao-no-intervalo}{\textit{item 2}} Proposition \ref{proposicao-propiedades-un})
\[u_{n,j}\rightarrow s_j \mbox{ in } C^{1,\alpha}[-j,j].\]

Let $\widetilde{a}_{j}:= \inf_{x\in [-j,j]} a_1(x)$. We will use the arguments presented on \textit{item 1} from Proposition \ref{proposicao-propiedades-un} to prove that $s_j$ \emph{is strictly positive on the interval }
$[-j,j]$.
\begin{center}
	\textbf{Case 1.} There exists a subsequence $(u_{n_i,j})_{i\in \bb{N}}$ of $(u_{n,j})$ such that $u'_{n_i,j}\geq 0$ in $(-j,0)$ for all $i$.
\end{center}
The analysis of this case follows exactly the same parameters of \textbf{Case 1} from \textit{Item 1}, Proposition \ref{proposicao-propiedades-un}. The main difference is the change of $\tilde{a}$ to $\widetilde{a}_{j}$.

\begin{center}
	\textbf{Case 2.}
	For all subsequence of $(u_{n,j})$ there exists a sub-subsequence $(u_{n_i,j})_{i\in \mathbb{N}}$ and exists a sequence $(z_i)_{i\in \mathbb{N}}\subset(-j,0)$ such that $u'_{n_i,j}(z_i)<0$.
\end{center}
For this case we can reformulate Lemma \ref{lema-caso2-limita-por-baixo} as follows: If $x\in (-j,j)$ and $u''_{n_i,j}(x)\geq 0$, then $u_{n_i,j}(x)> (\lambda\widetilde{a}_{j})^{\frac{1}{2-q}}$. This is true because we already know that $u_{n_i,j}$ is strictly positive in $(-j,j)$, then the estimation

\begin{equation}
	-A'(u_{n_i,j}(t))|u'_{n_i,j}(t)|^2-A(u_{n_i,j}(t))u''_{n_i,j}(t)+u_{n_i,j}(t)> \lambda\widetilde{a}_{j}|u_{n_i,j}(t)|^{q-1}
\end{equation}
is immediately established. The remaining argumentation is similar.

Thus, we conclude that $s_j>0$, as in Proposition \ref{proposicao-propiedades-un}. At last, let $\varphi \in C_0^{\infty}(-j,j)$. Then

\begin{multline*}
	\int_{-j}^{j}A(u_{n,j})u'_{n,j}\varphi' + \int_{-j}^{j}u_{n,j}\varphi = \int_{-j}^{j}(\lambda a_1 |u_{n,j}|^{q-1}+|u_{n,j}|^{p-1})\varphi \\ + \int_{-j}^{j} g(|u'_{n,j}|)\varphi.
\end{multline*}
When $n\rightarrow\infty$ we get
\begin{align*}
	\int_{-j}^{j}A(s_j)s_j'\varphi' + \int_{-j}^{j}s_j\varphi = \int_{-j}^{j}(\lambda a_1 |s_j|^{q-1}+|v|^{p-1})\varphi + \int_{-j}^{j} g(|s_j'|)\varphi.
\end{align*}
Since $\varphi\in C_0^{\infty}(-j,j)$ is arbitrary,  we conclude that  $s_j$ is a weak solution for the problem
\begin{equation}
	-(A(u)u')'(t) + u(t)= \lambda a_1(t) |u(t)|^{q-1} +  |u(t)|^{p-1} + g(|u'(t)|)
\end{equation}
in $(-j,j)$; by \cite[sect.~6.3, Theorem~1]{2010_Evans_BOOK} we have that $s_j\in H^2_{\text{loc}}(-j,j)$, thus -- using the same arguments as in (\ref{cara-da-segunda-derivada})-- $s_j\in C^2(-j,j)$.

\emph{Now we will construct our candidate solution of problem} (\ref{problema-principal-falandoSobre}). Define $w_n= u_{n,n}$, that is, $w_n$ is the diagonal sequence. Notice that, for $n\geq j$, $w_n$ is a subsequence of $u_{n,j}$; thus $w_n|_{[-j,j]}\rightarrow s_j$ in $C^{1,\alpha}[-j,j]$ for all $j\in \bb{N}$. Let $v(x)$ be defined by 
\[v(x)=\lim_{n\rightarrow \infty}w_n(x).\]
Then, $v(x)=s_j(x)$ for $x\in [-j,j]$. Since $\bb{R}=\bigcup_{j\in \bb{N}}[-j,j]$, by Lemma \ref{boa def das extensoes na reta}, $v$ is well defined in $\bb{R}$. Using the properties of $s_j$ obtained just above, and the fact that $v|_{[-j,j]}=s_j$ for all $j\in \bb{N}$, we conclude that $v\in C^2(\bb{R})$, $\|v\|_{H^1(\bb{R})}\leq \tilde{r}$ and $v$ is a positive, even solution of problem (\ref{problema-principal-falandoSobre}). From \cite[Pag.~214, Corol.~8.9]{2011_Brezis_BOOK} we get the homoclinic condition.

\subsection{Asymptotic solution }
Throughout our previous argumentation, we fixed $\lambda\in (0,\lambda^*]$. Denote by $v_\lambda$ the strong solution -- obtained above -- to the problem (\ref{problema-principal-falandoSobre}). We will analyze the behavior of $v_\lambda$ as $\lambda\rightarrow 0$.
\begin{prop}
	As $\lambda\rightarrow 0$, $v_\lambda \rightarrow 0$ in $C^0(\bb{R})$.
\end{prop}
\begin{proof}
    By \cite[Theorem 8.8.]{2011_Brezis_BOOK}, we obtain

 \begin{equation}\label{zc}
\|v_{\lambda}\|_{L^{\infty}(\mathbb{R})}\leq C \|v_{\lambda}\|_{H^1(\bb{R})}\leq C \tilde{r}.
 \end{equation}

	But remember that 
	\[\tilde{r} = \min\left\{ r, \lambda^{1/2}\left(\frac{2C_2 r^q}{\gamma}\right)^{1/2}   \right\}.\]
	Then, as $\lambda\rightarrow0 $, we obtain that $\|v_{\lambda}\|_{L^{\infty}(\mathbb{R})}\rightarrow 0$, completing the proof of Theorem \ref{teorema-central} .

\end{proof}


\subsection{Proof of Proposition \ref{prop1}}
We proceed to prove Proposition \ref{prop1} that says that there is no solution of \eqref{problema-principal-falandoSobre} for $\lambda$ large. 

\begin{proof}

 Suppose on the contrary that $\lambda^* = \infty$. In this way there is a sequence $\lambda_n \to \infty$ and corresponding solutions $v_{\lambda_n} > 0$ in $\mathbb{R}$ given by Theorem \ref{teorema-central}.

 Fix $R>0$ and define $\mathcal{P}(t,s) = \lambda a_1(t) s^{q-1}+s^{p-1}$ and $\widetilde{a}_R=\inf_{(-R,R)}a_1(t)$. Define also
$$
\Lambda = \lambda \  \widetilde{a}_R.
$$
We claim that there is a constant $C_{\Lambda} > 0$ such that
\begin{equation*}
\mathcal{P}(t,s) \geq \Lambda s^{q-1}+s^{p-1}\geq  C_{\Lambda} s \quad for \  \  s>0, \;\;  t\in(-R,R).
\end{equation*}

Consider the function $\mathcal{Q}(s) = (\Lambda s^{q-1}+s^{p-1}) s^{-1}$. Then $\mathcal{Q}(s) \to \infty$ as $s \to 0^+$ and as $s \to \infty$. The minimum value of $\mathcal{Q}$ is achieved at the unique point
$$
m = \left( \Lambda \  \frac{2 - q}{p - 2} \right)^\frac{1}{p - q}.
$$
Thus $C_{\Lambda} = \mathcal{Q}(m)$.

Let $\sigma_1>0$ and $\varphi_1>0$, respectively, the first eigenvalue and the first eigenfunction satisfying
$$
\left\{
\begin{array}{lll}
- \varphi_1'' = \sigma_1 \varphi_1 &\mbox{in}&(-R,R) \\
\varphi_1(-R)= \varphi_1(R)=0.
\end{array} \right.
$$
Since $C_{\Lambda}$ increases as $\lambda_n$ increases, there is $\lambda_0$  such that the corresponding constant satisfies $C_{\Lambda_0} \geq A(C\tilde{r})\sigma_1 + A(C\tilde{r})\delta + 1$, for all  $\delta \in(0,1)$. Hence, by $(H_2)-(H_3)$ and \eqref{zc}, the solution $v_{\lambda_0}>0$ in $\mathbb{R}$ of \eqref{problema-principal-falandoSobre} associated to $\lambda_0$ satisfies
$$
\left\{
\begin{array}{lll}
-  v''_{\lambda_0} \geq \frac{(C_{ \Lambda_0} - 1)}{A(v_{\lambda_0})} v_{\lambda_0}  \geq \frac{(C_{ \Lambda_0} - 1)}{A(C\tilde{r})} v_{\lambda_0} \geq (\sigma_1 + \delta) v_{\lambda_0} &\mbox{in}&(-R,R) \\
v_{\lambda_0}(-R), v_{\lambda_0}(R) \geq 0.
\end{array} \right.
$$
Otherwise, taking $\varepsilon>0$ small enough we obtain $\varepsilon \varphi_1 < v_0$ in $(-R,R)$ and
$$
\left\{
\begin{array}{lll}
- (\varepsilon \varphi_1)'' = (\varepsilon \sigma_1) \varphi_1 \leq (\sigma_1 + \delta) \varphi_1 &\mbox{in}&(-R,R) \\
\varphi_1(-R)=\varphi_1(R) =0.
\end{array} \right.
$$
By the method of subsolution and supersolution, there is a solution $\varepsilon \varphi_1 < \omega < v_0$ in $B_R(0)$ of
$$
\left\{
\begin{array}{lll}
- \omega'' = (\sigma_1 + \delta) \omega &\mbox{in}&(-R,R) \\
\omega(-R) =\omega(R) =0.
\end{array} \right.
$$
Hence there is a contradiction to the fact that $\sigma_1$ is isolated. Therefore, $\lambda^* < \infty$, indeed.
\end{proof}



  \bibliographystyle{elsarticle-num}
  \bibliography{referencias-edp}


%
%
%
\end{document}